\documentclass{amsart}
\usepackage{amssymb,amsmath, latexsym, enumerate, mathtools}
\usepackage[american,french]{babel}
\usepackage{amsopn}
\usepackage{graphicx}
\usepackage{fancyhdr}
\usepackage[T1]{fontenc}
\newtheorem{theorem}{Theorem}[section]

\newtheorem{lemma}[theorem]{Lemma}

\newtheorem{proposition}[theorem]{Proposition}

\def\fin{\hfill$\Box$\\}

\newcommand{\rr}{\mathbb{R}}
\newcommand{\eps}{\varepsilon}
\newcommand{\nn}{\mathbb{N}}
\newcommand{\cc}{\mathbb{C}}

\def\wrtext#1{\relax\ifmmode{\leavevmode\hbox{#1}}\else{#1}\fi}
\def\abs#1{\left|#1\right|}
\def\begeq{\begin{equation}}
\def\endeq{\end{equation}}

\def\part#1{\frac{\partial}{\partial #1}}

\newcommand{\real}{\mathbb{R}}

\renewcommand{\Re}{\textrm{Re }}
\renewcommand{\Im}{\textrm{Im }}

\def\og{orthogonal}

\def\wrt{with respect to }

\def\Re{{\rm Re\,}}
\def\Im{{\rm Im\,}}

\catcode`\é=\active\def é{\'e} \catcode`\è=\active\def è{\`e}
\catcode`\ç=\active\def ç{\c{c}} \catcode`\à=\active\def à{\`a}
\catcode`\ê=\active\def ê{\^e} \catcode`\ù=\active\def ù{\`u}
\catcode`\ô=\active\def ô{\^o} \catcode`\û=\active\def û{\^u}
\catcode`\î=\active\def î{\^{\i}} \catcode`\â=\active\def â{\^a}
\catcode`\ö=\active\def ö{\"o}

\def\mmm{{\mathcal M}} 

\def\sss{{\mathcal{S}}}

\def\R{\mathbb R}\def\N{\mathbb N}

\def\D{\partial}\def\eps{\varepsilon}\def\phi{\varphi}

\def\abs#1{\left\vert#1\right\vert}
\def\set#1{\left\{#1\right\}}\def\seq#1{\left<#1\right>}
\def\sep#1{\left(#1\right)}

\def\defegal{\stackrel{\text{\rm def}}{=}}

\begin{document}
\selectlanguage{american}
\title{ANISOTROPIC HYPOELLIPTIC ESTIMATES FOR LANDAU-TYPE OPERATORS}
\author{Fr\'ed\'eric H\'erau, Karel Pravda-Starov}

\address{\noindent \textsc{Laboratoire de Math\'ematiques,
UFR Sciences exactes et naturelles,
Universit\'e de Reims Champagne-Ardenne,
Moulin de la Housse,
BP 1039, 51687 Reims cedex 2, France}}
\email{herau@univ-reims.fr}
\urladdr{http://frederic.herau.perso.sfr.fr/}

\address{\noindent \textsc{
Department of Mathematics,
Imperial College London,
Huxley Building, 180 Queen's Gate,
London SW7 2AZ, UK}}
\email{k.pravda-starov@imperial.ac.uk}
\urladdr{http://www2.imperial.ac.uk/~kpravdas/index.html}

\begin{abstract}
We establish global hypoelliptic estimates for linear Landau-type operators. Linear Landau-type equations are a class of inhomogeneous kinetic equations with anisotropic diffusion whose study is motivated by the linearization of the Landau equation near the Maxwellian distribution.  By introducing a microlocal method by multiplier which can be adapted to various hypoelliptic kinetic equations, we establish for linear Landau-type operators optimal global hypoelliptic estimates with loss of $4/3$ derivatives in a Sobolev scale which is exactly related to the anisotropy of the diffusion.
\end{abstract}


\maketitle
\noindent

\section{Introduction}
An important problem in the theory of kinetic equations is concerned with studying the regularization properties of diffusive equations; and the derivation of sharp regularity estimates for their solutions. Among these equations are Fokker-Planck equations, Landau equations or Boltzmann equations without cut-off, either homogeneous or inhomogeneous.

Regarding the inhomogeneous case, that is, those kinetic equations describing the system evolution both in space and velocity variables, the analysis of these regularization phenomena is non-trivial, since diffusion generally occurs only in the velocity variable but not in the space one. In this sense, these equations can be considered as degenerate. Nevertheless, the regularization process in both space and velocity variables may still occur. This phenomenon essentially due to non-trivial interactions between the diffusive and transport parts of these equations, and known as hypoellipticity; is currently a very active domain of research in kinetic theory. We refer the reader to the series of recent works \cite{alex}, \cite{AMUXY09}, \cite{bouchut}, \cite{CLX09b}, \cite{HelNi05}, \cite{HerNi03}, \cite{MX07}, \cite{MX09}, \cite{xu}, which all highlight specific non-trivial mixing interactions between diffusion and transport leading to hypoellipticity in both space and velocity variables.

In the present work, we study the hypoellipticity of a particular class of inhomogeneous kinetic equations whose study is motivated by the linearization of the Landau equation near the Maxwellian distribution (see the end of this introduction).

We consider the class of linear Landau-type operators
\begin{equation}
P= iv.D_x+ D_v.\lambda(v)D_{v}+(v \wedge D_v).\mu(v) (v \wedge D_v)+F(v), \ x,v \in \rr^{3};\label{ind1}
\end{equation}
that is
\begin{multline*}
P = i\sum_{j=1}^3v_jD_{x_j}+ \sum_{j=1}^3D_{v_j}\lambda(v)D_{v_j}+(v_2D_{v_3}-v_3D_{v_2})\mu(v)(v_2D_{v_3}-v_3D_{v_2})\\
+(v_3D_{v_1}-v_1D_{v_3})\mu(v)(v_3D_{v_1}-v_1D_{v_3}) \\ +(v_1D_{v_2}-v_2D_{v_1})\mu(v)(v_1D_{v_2}-v_2D_{v_1})+F(v),
\end{multline*}
with $D_x=i^{-1}\partial_x$, $D_v=i^{-1}\partial_v$ and $\gamma \in [-3,1]$;  where the diffusion is given by smooth positive functions $\lambda$, $\mu$ and $F$ satisfying for all $\alpha \in \nn^3$,
\begin{equation}\label{eq-3.1}
\exists C_{\alpha}>0, \forall v \in \rr^3,  |\partial_v^{\alpha}\lambda(v)| +|\partial_v^{\alpha}\mu(v)| \leq C_{\alpha} \langle v \rangle^{\gamma-|\alpha|}; |\partial_v^{\alpha}F(v)| \leq C_{\alpha} \langle v \rangle^{\gamma+2-|\alpha|};
\end{equation}
and
\begin{equation}\label{w2}
\exists C>0, \forall v \in \rr^3,  \ \lambda(v) \geq C \langle v \rangle^{\gamma}; \ \mu(v) \geq C \langle v \rangle^{\gamma}; \  F(v) \geq C \langle v \rangle^{\gamma+2};
\end{equation}
with $\langle v \rangle=(1+|v|^2)^{\frac{1}{2}}$. Linear Landau-type operators are formally accretive operators
$$\textrm{Re}(Pu,u)_{L^2}=\|\lambda(v)^{\frac{1}{2}}D_vu\|_{L^2}^2+\|\mu(v)^{\frac{1}{2}}(v \wedge D_v)u\|_{L^2}^2+\|F(v)^{\frac{1}{2}}u\|_{L^2}^2 \geq 0, \ u \in \mathcal{S}(\rr_{x,v}^6);$$
with an anisotropic diffusion due to the presence of the cross product term $v \wedge D_v$. Denoting $(\xi,\eta)$ the dual variables of $(x,v)$, we notice that the diffusion only occurs in the variables $(v,\eta)$, but not in the other directions; and that the cross product term $v \wedge D_v$ improves this diffusion in specific directions of the phase space where the variables $v$ and $\eta$ are orthogonal. In this work, we aim at proving that linear Landau-type operators are actually hypoelliptic despite this lack of diffusion in the spatial derivative~$D_x$. More specifically, we shall be concerned in proving optimal global hypoelliptic estimates in a specific Sobolev scale in both spatial and velocity derivatives whose structure is exactly related to the anisotropy of the diffusion.

The main result of this article is given by the following  global anisotropic hypoelliptic estimate with loss of $4/3$ derivatives:

\bigskip

\begin{theorem} \label{mainlandau} Let $P$ be the linear Landau-type operator defined in \emph{(\ref{ind1})}. Then, there exists a positive constant $C > 0$ such that for all $u \in \sss(\R^{6}_{x,v})$,
\begin{multline}\label{ind2}
\|\langle v \rangle^{\gamma +2} u\|_{L^2}^2 +  \|\langle v \rangle^{\gamma} |D_v|^{2} u \|_{L^2}^2 +  \| \langle v \rangle^{\gamma} |v \wedge D_v|^{2} u \|_{L^2}^2 \\
+ \|\langle v\rangle^{\gamma/3} |D_x|^{2/3} u \|_{L^2}^2 +  \|\langle v\rangle^{\gamma/3} |v \wedge D_x|^{2/3} u \|_{L^2}^2
  \leq C  (\|Pu\|_{L^2}^2 + \|u\|_{L^2}^2),
\end{multline}
where the notation $\|\cdot\|_{L^2}$ stands for the $L^2(\rr_{x,v}^{6})$-norm.
\end{theorem}

\bigskip

We begin by noticing that the terms controlled in this global hypoelliptic estimate are sharp and have an anisotropic structure similar to the diffusion term. More specifically, as in the diffusion term, the presence of the two cross products $v \wedge D_v$ and $v \wedge D_x$ in
$$\| \langle v \rangle^{\gamma} |v \wedge D_v|^{2} u \|_{L^2}^2+ \|\langle v\rangle^{\gamma/3} |v \wedge D_x|^{2/3} u \|_{L^2}^2,$$
improves the regularity estimates provided by the terms
$$\|\langle v \rangle^{\gamma} |D_v|^{2} u \|_{L^2}^2+\|\langle v\rangle^{\gamma/3} |D_x|^{2/3} u \|_{L^2}^2,$$
in the specific directions of the phase space where either, $v$ and $D_v$, or $v$ and $D_x$; are orthogonal. The anisotropy and the different indices appearing in the estimate (\ref{ind2}) are optimal. Notice indeed that this hypoelliptic estimate splits up into two parts. The first part of the estimate
\begin{equation}\label{ind3}
\|\langle v \rangle^{\gamma +2} u\|_{L^2}^2 +  \|\langle v \rangle^{\gamma} |D_v|^{2} u \|_{L^2}^2 +  \| \langle v \rangle^{\gamma} |v \wedge D_v|^{2} u \|_{L^2}^2
  \leq C  (\|Pu\|_{L^2}^2 + \|u\|_{L^2}^2),
\end{equation}
is purely provided by the diffusion term of the linear Landau-type operator; and we notice from (\ref{w2}) that the left-hand-side of (\ref{ind3}) has exactly the same anisotropic structure and asymptotic growth as the diffusion term
$$D_v.\lambda(v)D_{v}+(v \wedge D_v).\mu(v) (v \wedge D_v)+F(v).$$
It follows that this first part of the estimate is obviously optimal. On the other hand, the most interesting result in Theorem~\ref{mainlandau} is the anisotropic regularity estimate in the spatial derivative $D_x$,
\begin{equation}\label{ind4}
\|\langle v\rangle^{\gamma/3} |D_x|^{2/3} u \|_{L^2}^2 +  \|\langle v\rangle^{\gamma/3} |v \wedge D_x|^{2/3} u \|_{L^2}^2
  \leq C  (\|Pu\|_{L^2}^2 + \|u\|_{L^2}^2).
\end{equation}
This second estimate is also optimal in term of the index $2/3$ appearing in the left-hand-side of (\ref{ind4}). Indeed, the optimality of this index $2/3$ is suggested by general results about microlocal hypoellipticity with optimal loss of derivatives established in~\cite{camus} (Corollary~1.3) or \cite{phong}. Let us recall the general result about microlocal hypoellipticity proved by P.~Bolley, J.~Camus and J.~Nourrigat in~\cite{camus} (Theorem~1.1 and Corollary~1.3):
Let $(A_j)_{1\leq j\leq l}$ be a system of properly supported classical pseudodifferential operators $(\rho=1,\delta=0)$ on an open subset $\Omega$ of $\rr^n$ of arbitrary real orders $m_1,\cdots,m_l$. Suppose that $A_j-A_j^\ast$ has order $m_j-1$ for all $1 \leq j \leq l$. Let $(x_0,\xi_0) \in T^\ast(\Omega)\smallsetminus\{0\}$ be such that there is a commutator of length~$r$, $Y=(\text{ad}\,A_{i_1})\cdots(\text{ad}\,A_{i_{r-1}})A_{i_r}$, which is elliptic of order $m_{i_1}+\cdots+m_{i_r}-r+1$ at $(x_0,\xi_0)$. Then the following implication holds for all $s\in \rr$: If $u\in \mathcal{D}'(\Omega)$ and $A_ju\in H^{s-m_j}(x_0,\xi_0)$, $j=1,\cdots,l$; then $u\in H^{s-1+1/r}(x_0,\xi_0)$. As a corollary, one obtains that if all the $m_j$ are equal then $\Sigma_1^lA_j^\ast A_j$ is hypoelliptic at $(x_0,\xi_0)$ with loss of $2(1-r^{-1})$ derivatives.
When each $A_j$ is a real vector field, this is a microlocal version of the celebrated theorem by L.~H\"ormander on the hypoellipticity of \og sums of squares\fg \ proved in~\cite{Hor67}. A simpler proof of H\"ormander Theorem, but with less precise information on the loss of derivatives, was given by J.J.~Kohn in~\cite{kohn}; whereas optimal estimates for the loss of derivatives were obtained, in the case of real vector fields, by L.P.~Rothschild and E.M.~Stein in~\cite{stein}. Linear Landau-type operators are non-selfadjoint operators for which these general results of hypoellipticity does not apply. However, as mentioned above, hypoellipticity for linear Landau-type operators will be derived from non-trivial mixing interactions between their diffusion and transport parts. More specifically, hypoellipticity for linear Landau-type operators will come from the ellipticity of commutators of length 3 of their diffusion and transport parts. This explains that the optimal loss of derivatives expected in this case is $2(1-1/3)=4/3$; and that the order $2$ associated to the diffusion term and the regularity estimate with respect to the velocity derivative $D_v$ must be substituted by an order $2-4/3=2/3$ in the regularity estimate with respect to the spatial derivative $D_x$. Regarding now the anisotropic structure of the term appearing in the left-hand-side of the estimate (\ref{ind4}), this structure will directly come from the explicit expression of the Poisson brackets associated to these elliptic commutators of length~3.

Kohn's method is the simplest and most flexible way for proving hypoellipticity. However, it does not provide the optimal loss of derivatives. In order to obtain the optimal loss of derivatives, more subtle microlocal and geometric methods are needed. In this work, we shall present a general method by multiplier which allows to prove hypoellipticity with optimal loss of $4/3$ derivatives. This method has been first introduced by F.~H\'erau, J.~Sj\"ostrand and C.~Stolk in their work on the Fokker-Planck equation~\cite{HSS06}. This approach has then been extended in a specific case~\cite{Pra09} by the second author to get optimal hypoelliptic estimates with loss of $2(1-(2k+1)^{-1})$, $k \in \nn$, derivatives. Because this method is very general and that it can be adapted to various hypoelliptic kinetic equations, we aim here at giving an extensive presentation of this approach. In order to do so, we shall first apply this method (Section~\ref{fok}) to recover the well-known hypoellipticity with loss of $4/3$ derivatives for the Fokker-Planck operator without external potential
$$P=iv.D_x+D_v^2+v^2.$$
This example of the Fokker-Planck operator will allow to present the principles of this multiplier method in a simplified setting where there is a good symbolic calculus.
In a second step, we shall then consider linear Landau-type operators and prove Theorem~\ref{mainlandau} (Section~\ref{la1}). We will see that this general multiplier method is sharp enough to handle anisotropic classes of symbols. However, because of this anisotropy, we will have to deal with gainless symbolic calculus. As a consequence, the implementation of this method in the case of linear Landau-type operators will be more complex and will require the use of more advanced microlocal analysis. In order to handle this setting with gainless symbolic calculus, we shall use some elements of Wick calculus developed by N.~Lerner in~\cite{Ler03}. For convenience of reading, the main features and the definition of Wick calculus is recalled in a short self-contained presentation given in appendix (Section~\ref{la3}).

Finally, we shall end this introduction by giving few elements of explanations about the motivation for studying this class of linear Landau-type operators.
Linear Landau-type equations are a class of inhomogeneous kinetic equations whose study is motivated by the linearization of the Landau equation. Details about the Landau equation may be found for example in the works by Y.~Guo~\cite{Guo02}, C.~Mouhot and L.~Neumann~\cite{MN06}, or C.~Villani~\cite{Vil02}; and we may only recall here that the Landau equation reads as the evolution equation of the density of particles
\begin{equation} \label{QL}
\begin{cases}
\partial_tf+v\cdot\nabla_{x}f=Q_L(f,f),\\
f|_{t=0}=f_0
\end{cases}
\end{equation}
where $Q_L$ is the so-called Landau collision operator
\begin{equation} \label{land}
Q_L(f,f) = \nabla_v \cdot \Big(\int_{\R^3} {\bf A}(v-v_*)\big(f(v_*) (\nabla_v f)(v) - f(v) (\nabla_vf)(v_*)\big) dv_*\Big).
\end{equation}
 Here, $\textbf{A}(z)$ is a symmetric nonnegative matrix depending on a parameter $z \in \R^3$,
$${\bf A}(z) = |z|^2 \Phi(|z|) {\bf P}(z),$$
 with   $ \Phi(|z|) = |z|^\gamma$  and $\gamma \in [-3,1]$;
 which is proportional to $\textbf{P}$ the orthogonal projection onto $z^{\perp}$,
 $$\textbf{P}(z) = \textrm{Id} - \frac{1}{|z|^2} z. z^\perp,$$
 matrix whose entries are
 $$\big(\textbf{P}(z)\big)_{i,j} = \delta_{i,j} - \frac{z_i z_j}{|z|^2}, \ 1 \leq i,j \leq 3.$$
 The original Landau collision operator describing collisions among charged particles interacting with Coulombic force and introduced by Landau in~1936, corresponds to the case $\gamma=-3$. As in the Boltzmann equation, it is well-known that Maxwellians are steady states to the Landau equation
 \begin{equation} \label{max}
 \mmm(x,v)=(2 \pi)^{-3/2}e^{-|v|^2/2}.
\end{equation}
Following the standard procedure described in~\cite{Guo02} or \cite{MN06},  we linearize the Landau equation around $\mmm$
 by posing
$$f=\mmm+\sqrt{\mmm}u,$$
and one can check that  \it after linearization \rm the Landau equation for the perturbation $u(t,x,v)$ now reads as
 \begin{equation} \label{le}
\D_t u + iv. D_xu -L u = 0,
\end{equation}
with $D_x=i^{-1}\partial_x$. The transport part of the equation $iv.D_x$ is unchanged, whereas one can prove that the operator $L$ may write as
 \begin{equation} \label{landaulin0}
 L= L_* - D_v A (v) D_v - F(v),
\end{equation}
with $F$ a positive smooth function satisfying the estimates (\ref{eq-3.1}) and (\ref{w2}). Here, the operator $L_*$ is a convolution-type term bounded on $L^2$, which only has a (big) influence on the lower part of the spectrum of the operator $iv.D_x-L;$
whereas the other term
\begin{equation} \label{landaulin}
A (v) = ({\bf A} * \mmm)(v),
\end{equation}
inherits the properties of the projection ${\bf P}$. More specifically, for each vector $v \in \rr^3$, the matrix $A(v)$ is symmetric
with a simple eigenvalue  $\lambda(v)$ associated to the eigenvector $v$; and a double eigenvalue $\lambda_\perp(v)$ associated to the eigenspace $v^\perp$; which satisfy the estimates
$$\forall \alpha \in \nn^3, \exists C_{\alpha}>0, \forall v \in \rr^3, \ |\partial_v^{\alpha} \lambda(v)| \leq C_\alpha \langle v\rangle^{\gamma-|\alpha|}; \ |\partial_v^{\alpha} \lambda_\perp(v)|  \leq C_\alpha \langle v \rangle^{\gamma+2-|\alpha|},
 $$
giving rise to the anisotropy of the diffusion. Up to a bounded operator, this explains why the linearization of the Landau equation essentially reduces to the study of a linear Landau-type operator
$$P = iv.D_x+D_v.\lambda(v)D_{v}+(v \wedge D_v).\mu(v)(v \wedge D_v)+ F(v),$$
with
$\mu(v) \sim  \frac{\lambda_\perp(v)}{\seq{v}^2};$ and a perhaps slightly modified function $\lambda(v)$
so that the estimates (\ref{w2}) hold. This motivates the present work on the hypoellipticity of these operators.

\section{Optimal hypoelliptic estimate for the Fokker-Planck operator}\label{fok}

As mentioned in the introduction, we shall first consider the case of the Fokker-Planck operator without external potential
\begin{equation}\label{ind10}
P=iv.D_x+D_v^2+v^2, \ x,v \in \rr^n;
\end{equation}
which provides a neat setting for explaining the principles of the general method we shall use later on for proving the hypoellipticity of linear Landau-type operators.

More specifically, we aim in this section at recovering the following well-known optimal hypoelliptic estimate with loss of $4/3$ derivatives:

\bigskip

 \begin{proposition}\label{KFP}
 Let $P$ be the Fokker-Planck operator defined in \emph{(\ref{ind10})}. Then, there exists a positive constant $C>0$ such that for all $u \in \sss(\R^{2n}_{x,v})$,
 $$\|\langle D_x \rangle^{2/3}u\|_{L^2}^2 +\|\langle v \rangle^{2} u\|_{L^2}^2 + \|\langle D_v \rangle^{2}u\|_{L^2}^2 \leq C(\|Pu\|_{L^2}^2 + \|u\|_{L^2}^2),$$
where the notation $\|\cdot\|_{L^2}$ stands for the $L^2(\rr_{x,v}^{2n})$-norm.
\end{proposition}

\bigskip

\noindent
This result of hypoellipticity is essentially contained in~\cite{HSS06} (Sections 2, 8 and 9); and we shall use this example of the Fokker-Planck operator as a model to illustrate in a simplified setting with good symbolic calculus a general method for proving optimal hypoelliptic estimates with loss of $4/3$ derivatives. This microlocal method by multiplier can be adapted to various hypoelliptic kinetic equations; and as we shall see with linear Landau-type operators, it turns out to be sharp enough to handle anisotropic classes of symbols, even if in the latter case we shall have to deal with gainless symbolic calculus.

Coming back from now to the Fokker-Planck operator, we begin by performing a partial Fourier transform in the $x$ variable; and notice that one may reduce our study on the Fourier side to the analysis of the operator
$$P=iv.\xi + D_v^2 + v^2=iv.\xi+ \sum_{j=1}^nD_{v_j}^2+ \sum_{j=1}^n v_j^2, \  v,\xi \in \rr^{n};$$
depending on the parameter $\xi$. In this section, we shall therefore consider Weyl quantizations of symbols only in the velocity variable $v$ and its dual variable $\eta$; but not in the variable $\xi$, which will be considered here as a parameter
\begin{equation}\label{fili2}
(a^wu)(v)=\frac{1}{(2\pi)^n}\int_{\rr^{2n}}e^{i(v-\tilde{v}).\eta}a\Big(\frac{v+\tilde{v}}{2},\eta\Big)u(\tilde{v})d\tilde{v}d\eta.
\end{equation}
The Weyl symbol of the Fokker-Planck operator is then given by
$$p=iv.\xi+|\eta|^2+|v|^2,$$
where $|\cdot|$ stands for the Euclidean norm on $\rr^n$.
Defining the symbol
\begin{equation}\label{m2fp}
\lambda=\big(1+|\eta|^2+|v|^2+|\xi|^2\big)^{\frac{1}{2}},
\end{equation}
we shall see that Proposition~\ref{KFP} easily follows from the key hypoelliptic estimate
\begin{equation}\label{ya1}
\|(\lambda^{2/3})^wu\|_{L^2}^2 \lesssim \|Pu\|_{L^2}^2+\|u\|_{L^2}^2.
\end{equation}
In order to explain how one can derive such an hypoelliptic estimate and justify the choice of multiplier introduced below, we first notice that the diffusive part of the Fokker-Planck operator  gives a trivial control in the variables $(v,\eta)$. Indeed, this control is just a consequence of the ellipticity of the real part of the symbol
$$\textrm{Re }p=|\eta|^2+|v|^2,$$
in these variables.
The main point in the estimate (\ref{ya1}) is then to get a control of the term $|\xi|^{2/3}$. Notice that this control cannot be derived from the ellipticity of the symbol $p$; and that we will need to consider the following iterated commutator
$$[(\textrm{Im }p)^w,[(\textrm{Re }p)^w,(\textrm{Im }p)^w]],$$
where $\textrm{Re }p$ and $\textrm{Im }p$ stand for the real and imaginary parts of the symbol $p$; in order to get some ellipticity in the parameter $\xi$.
Indeed, usual symbolic calculus (see Theorem~18.5.4 in~\cite{Hor85}) or a direct computation shows that the Weyl symbol of this iterated commutator is exactly given by the iterated Poisson brackets
$$-\{\textrm{Im }p,\{\textrm{Re }p,\textrm{Im }p\}\}=\{\textrm{Im }p,\{\textrm{Im }p,\textrm{Re }p\}\}=2|\xi|^2.$$
The Poisson bracket of two symbols $a$ and $b$ is defined as
$$\{a,b\}=H_ab=\frac{\partial a}{\partial \eta}\cdot \frac{\partial b}{\partial v}-\frac{\partial a}{\partial v}\cdot \frac{\partial b}{\partial \eta},$$
where $H_a$ stands for the Hamilton vector field of $a$,
$$H_a=\frac{\partial a}{\partial \eta}\cdot \frac{\partial }{\partial v}-\frac{\partial a}{\partial v}\cdot \frac{\partial }{\partial \eta}.$$
Notice that we shall need the ellipticity of this iterated commutator only in the region of the phase space where $|\eta|^2+|v|^2 \lesssim \lambda^{2/3}$;
since one can directly rely on the real part of the symbol $p$ in the region where $|\eta|^2+|v|^2 \gtrsim \lambda^{2/3}$. This informal discussion accounts for the following choice of symbol multiplier.
Let $\psi$ be a $C_0^{\infty}(\rr,[0,1])$ function such that
\begin{equation}\label{n5fp}
\psi=1 \textrm{ on } [-1,1], \textrm{ and } \textrm{supp }\psi \subset [-2,2].
\end{equation}
We define the real-valued symbol
\begin{equation}\label{m1fp}
g=-\frac{\xi.\eta}{\lambda^{4/3}}
\psi\left(\frac{|\eta|^2+|v|^2}{\lambda^{2/3}}\right),
\end{equation}
where the function $\lambda$ is defined in (\ref{m2fp}). The cutoff function $\psi$ allows to localize the symbol multiplier in the region of the phase space where we need the ellipticity of the iterated commutator
$$[(\textrm{Im }p)^w,[(\textrm{Re }p)^w,(\textrm{Im }p)^w]].$$
It is essential to localize the symbol multiplier exactly in this region if we want to get the optimal loss of derivatives in the hypoelliptic estimate (\ref{ya1}). Notice that the term $\xi.\eta$ appearing in the expression of the symbol $g$ will play an essential r\^ole in the following. Up to a factor 2, it is actually equal to the symbol
$$H_{\textrm{Re}p}\textrm{ Im }p=\{\textrm{Re }p,\textrm{Im }p\}=2 \xi.\eta.$$
As we shall see below, this term will make appear the elliptic symbol of the iterated commutator
$$-H_{\textrm{Im}p}H_{\textrm{Re}p} \textrm{ Im }p=H_{\textrm{Im}p}^2 \textrm{ Re }p=2|\xi|^2;$$
whereas the factor $\lambda^{4/3}$ appearing in (\ref{m1fp}) will ensure that the symbol $g$ defines a bounded operator on $L^2$.
Following the usual notations introduced by L.~H\"ormander in~\cite{Hor85} (Chapter~18), see also \cite{Le}; we consider the metric
$$\Gamma = \frac{dv^2+d\eta^2}{M},$$
with
\begin{equation}\label{ya2}
M = 1+|v|^2 + |\eta|^2 + \lambda^{2/3};
\end{equation}
and the classes of symbols $S(m,\Gamma)$ associated to order functions $m$, that is, the class of all functions $a\in C^{\infty}(\real_{v,\eta}^{2n},\cc)$ possibly depending on the parameter $\xi$; and satisfying
$$\forall \alpha \in \nn^{2n}, \exists C_{\alpha}>0, \forall (v,\eta,\xi) \in \rr^{3n}, \  |\partial_{v,\eta}^{\alpha} a(v,\eta,\xi)| \leq C_{\alpha} m(v,\eta,\xi) M(v,\eta,\xi)^{-|\alpha|/2}.$$
It is easy to check that this metric $\Gamma$ is admissible (slowly varying, satisfying the uncertainty principle and temperate) with gain
\begin{equation}\label{dee1}
\lambda_{\Gamma}(X)=\inf_{T \neq 0}\left(\frac{\Gamma_X^{\sigma}(T)}{\Gamma_X(T)}\right)^{1/2}=M(X), \ X=(v,\eta,\xi);
\end{equation}
for symbolic calculus in the symbol classes $S(m,\Gamma)$. We refer to \cite{Hor85} or \cite{Le} for extensive presentations of symbolic calculus. We begin by proving the following symbolic estimates:

\bigskip

\begin{lemma} \label{2.2fp}
For any $m \in \R$, the following symbols belong to their respective symbol classes
\begin{align*}
i) & \ \langle \xi \rangle^m \in S( \lambda^m, \Gamma); \ \ ii)\ \lambda^m \in S( \lambda^m, \Gamma); \ \ iii)\ \psi\left(\frac{|\eta|^2+|v|^2}{\lambda^{2/3}}\right) \in S(1,\Gamma);\\
iv) & \  g \in S(1,\Gamma); \ \ v) \ \Re p  \in S(M,\Gamma);
\end{align*}
uniformly \wrt the parameter $\xi \in \R^n$.
\end{lemma}

\bigskip

\proof The assertion $i)$ is obvious since the term $\langle \xi \rangle^m$ is independent of the variables $(v,\eta)$. By using the writing convention
$$f(X) \lesssim g(X),$$
with $X \in \rr^d$; for the existence of a positive function $C>0$ such that the estimate
$$f(X) \leq C g(X),$$
holds for all $X \in \rr^d$; we easily notice from (\ref{m2fp}) and (\ref{ya2}) that for all $\alpha  \in \N^{2n}$, we have
$$|\D_{v,\eta}^\alpha(\lambda^m)| \lesssim \lambda^{m -|\alpha|} \lesssim \lambda^m M^{-|\alpha|/2},$$
uniformly \wrt the parameter $\xi \in \R^n$; since the estimate $M^{1/2} \lesssim \lambda$ holds uniformly \wrt $\xi$. This proves assertion $ii)$. Regarding assertion $iii)$, we first notice that on the support of the function
$$\psi\left(\frac{|\eta|^2+|v|^2}{\lambda^{2/3}}\right),$$
the estimate $|\eta|^2 + |v|^2 \lesssim \lambda^{2/3}$ implies that
$$M^{1/2} \sim \lambda^{1/3}$$
and
$$|\D^\alpha_{v,\eta} (|\eta|^2 + |v|^2)|  \lesssim \left\{ \begin{array}{ll} \lambda^{2/3} & \textrm{ when } |\alpha| = 0,\\  \lambda^{1/3} & \textrm{ when } |\alpha| = 1,\\ 1 & \textrm{ when } |\alpha| = 2, \\ 0 & \textrm{ when } |\alpha| \geq 3. \end{array} \right.$$
Assertion $iii)$ then directly follows assertion from $ii)$. We next notice that on the support of the function
$$\psi\left(\frac{|\eta|^2+|v|^2}{\lambda^{2/3}}\right),$$
the estimate $|\xi.\eta| \leq |\xi||\eta| \lesssim \lambda^{4/3}$ implies that
\begin{equation}\label{fili1}
|\D^\alpha_{v,\eta} (\xi.\eta)|  \lesssim \left\{ \begin{array}{ll} \lambda^{4/3} & \textrm{ when } |\alpha| = 0,\\  \lambda & \textrm{ when } |\alpha| = 1,\\ 0 & \textrm{ when } |\alpha| \geq 2. \end{array} \right.
\end{equation}
Recalling that in this region $M^{1/2} \sim \lambda^{1/3}$, assertion $iv)$ is then a direct consequence of assertions $ii)$ and $iii)$; whereas assertion $v)$ is trivial.  \fin

\noindent
Next Lemma shows that up to controlled terms and a weight factor $\lambda^{4/3}$, the Poisson bracket
$$H_{\textrm{Im}p}\ g=\set{ \Im p, g},$$
makes appear the elliptic symbol of the iterated commutator
$$-H_{\textrm{Im}p}H_{\textrm{Re}p} \textrm{ Im }p=H_{\textrm{Im}p}^2 \textrm{ Re }p=2|\xi|^2,$$
in the region of the phase space where $|\eta|^2+|v|^2 \lesssim \lambda^{2/3}$.

\bigskip

\begin{lemma} \label{mainterm}
We have
$$
H_{\emph{\textrm{Im}}p}\ g = \frac{|\xi|^2}{\lambda^{4/3}}
\psi\left(\frac{|\eta|^2+|v|^2}{\lambda^{2/3}}\right) + r,
$$
with a remainder $r$ belonging to both symbol classes $S\big(|\eta|^2+|v|^{2},\Gamma\big)$ and $S(M,\Gamma)$, uniformly \wrt the parameter $\xi \in \R^n$.
\end{lemma}

\bigskip

\proof Recalling the definition (\ref{m1fp}), an explicit computation of the Poisson bracket
$$H_{\textrm{Im}p}\ g =\{\textrm{Im }p,g\}=\{\xi.v,g\},$$
gives that
\begin{equation}\label{n7fp}
\{\xi.v,g\}=-\xi.\frac{\partial g}{\partial \eta}=
\frac{|\xi|^2}{\lambda^{4/3}}\psi\left(\frac{|\eta|^2+|v|^2}{\lambda^{2/3}}\right)+r
\end{equation}
with
\begin{equation}
r=\big(\xi.\eta\big)\big(\xi.\partial_{\eta}\big(\lambda^{-4/3}\big)\big)\psi\left(\frac{|\eta|^2+|v|^2}{\lambda^{2/3}}\right)
+\frac{\xi.\eta}{\lambda^{4/3}}(\xi.\partial_{\eta})\left[\psi\left(\frac{|\eta|^2+ |v|^2}{\lambda^{2/3}}\right)\right].
\end{equation}
Recalling that $M^{1/2} \sim \lambda^{1/3}$ on the support of the function
$$\psi\left(\frac{|\eta|^2+|v|^2}{\lambda^{2/3}}\right),$$
we then notice from Lemma~\ref{2.2fp} and (\ref{fili1}) that the term
\begin{align*}
\big(\xi.\eta\big)\big(\xi.\partial_{\eta}\big(\lambda^{-4/3}\big)\big) \psi\left(\frac{|\eta|^2+|v|^2}{\lambda^{2/3}}\right) & = -\frac{4}{3} \big(\xi.\eta\big)^2 \lambda^{-10/3} \psi\left(\frac{|\eta|^2+|v|^2}{\lambda^{2/3}}\right) \\
& \in S(\lambda^{-2/3},\Gamma)
\end{align*}
and
$$\frac{\xi.\eta}{\lambda^{4/3}}(\xi.\partial_{\eta})\left[\psi\left(\frac{|\eta|^2+ |v|^2}{\lambda^{2/3}}\right)\right] \in S(\lambda^{2/3},\Gamma),$$
since $\xi \in S(\lambda,\Gamma)$. By using now that
$$|\eta|^2+|v|^2 \sim \lambda^{2/3},$$
on the support of the function
$$\psi'\left(\frac{|\eta|^2+|v|^2}{\lambda^{2/3}}\right),$$
we deduce that
$$\big(\xi.\eta\big)\big(\xi.\partial_{\eta}\big(\lambda^{-4/3}\big)\big) \psi\left(\frac{|\eta|^2+|v|^2}{\lambda^{2/3}}\right) \in S(1,\Gamma)$$
and
$$\frac{\xi.\eta}{\lambda^{4/3}}(\xi.\partial_{\eta})\left[\psi\left(\frac{|\eta|^2+ |v|^2}{\lambda^{2/3}}\right)\right] \in S\big(|\eta|^2+|v|^{2},\Gamma\big),$$
uniformly \wrt the parameter $\xi \in \R^n$. Recalling that $|\eta|^2+|v|^{2} \leq M$, we finally obtain that the remainder $r$ belongs to both symbol classes $S\big(|\eta|^2+|v|^{2},\Gamma\big)$ and $S(M,\Gamma)$, uniformly \wrt the parameter $\xi \in \R^n$.
\fin

\noindent
By using Lemma~\ref{mainterm}, we can then prove the following estimate:

\bigskip

 \begin{proposition} \label{KFP2} There exists a positive constant $C>0$ such that for all $s \in \rr$, $\xi \in \rr^n$ and $u \in \sss (\R^{n}_{v})$,
 $$ \||\xi|^{1/3}u\|_{L^2}^2 +\|v u\|_{L^2}^2+ \|D_v u\|_{L^2}^2 \leq C \big( \|\langle \xi\rangle^{-s}Pu \|_{L^2} \|\langle\xi\rangle^s u\|_{L^2}+\|u\|_{L^2}^2\big),$$
 where $\|\cdot\|_{L^2}$ stands for the $L^2(\rr_{v}^{n})$-norm.
\end{proposition}

\bigskip

\proof We consider the multiplier $G = g^w$ defined by the Weyl quantization of the symbol $g$ as in (\ref{fili2}); and let $\eps$ be a positive parameter such that $0<\eps \leq 1$. For any $s \in \rr$ and $0<\eps \leq 1$, we may write
\begin{multline}\label{eq7fp}
\textrm{Re}(\langle  \xi \rangle^{-s}Pu,\langle  \xi \rangle^{s}(1-\eps G)u)=\|D_vu\|_{L^2}^2+\|v u\|_{L^2}^2\\ -\eps \textrm{Re}(iv.\xi u,Gu)
-\eps \textrm{Re}(|D_v|^2 u,Gu)-\eps \textrm{Re}(|v|^2 u,Gu).
\end{multline}
 We need to estimate the terms appearing on the second line of (\ref{eq7fp}). We begin by noticing from Lemma~\ref{2.2fp} and the Calder\'on-Vaillancourt Theorem that the operator $G$ is bounded on $L^2$. This implies that
\begin{multline}
|\textrm{Re}(|D_v|^2 u,Gu)|=|\textrm{Re}(D_v u,D_v Gu)| \\ \leq |\textrm{Re}(D_v u,[D_v,G]u)|+|\textrm{Re}(D_v u,G D_vu)|
\lesssim \|D_v u\|_{L^2}^2+\|[D_v,G]u\|_{L^2}^2,
\end{multline}
uniformly \wrt the parameter $\xi \in \R^n$.
Symbolic calculus shows that the symbol of the commutator $[D_v,G]$ is exactly given by $i^{-1}\partial_v{g}$. In view of Lemma~\ref{2.2fp}, this symbol
belongs to the symbol class $S(1,\Gamma)$.  We therefore deduce from the Calder\'on-Vaillancourt Theorem that
\begin{equation}\label{fili3}
|\textrm{Re}(|D_v|^2 u,Gu)| \lesssim \|D_v u\|_{L^2}^2+\|u\|_{L^2}^2,
\end{equation}
uniformly \wrt the parameter $\xi \in \R^n$.
A similar reasoning gives the estimate
\begin{equation}\label{fili4}
|\textrm{Re}(|v|^2 u,Gu)| \lesssim \|vu\|_{L^2}^2+\|u\|_{L^2}^2,
\end{equation}
uniformly \wrt the parameter $\xi \in \R^n$.
Regarding the last term, we may write
$$-\textrm{Re}(iv.\xi u,Gu)= \frac{1}{2}\textrm{Re}([i v.\xi ,G]u,u),$$
since the operators $G$ and $iv.\xi$ are respectively formally selfadjoint and skew-selfadjoint.
Symbolic calculus then shows that the symbol of the commutator
$$\frac{1}{2}[i v.\xi ,G],$$
is exactly given by
$$\frac{1}{2}H_{\textrm{Im}p}\ g=\frac{1}{2}\set{v.\xi, g}.$$
Lemma~\ref{mainterm} shows that the symbol of this commutator may be written as
$$H_{\textrm{Im}p}\ g = \frac{|\xi|^2}{\lambda^{4/3}}\psi\left(\frac{|\eta|^2+|v|^2}{\lambda^{2/3}}\right) + r,$$
where $r$ stands for a remainder belonging to both symbol classes $S\big(|\eta|^2+|v|^{2},\Gamma\big)$ and $S(M,\Gamma)$, uniformly \wrt the parameter $\xi \in \R^n$.
Notice from Lemma~\ref{2.2fp} and (\ref{dee1}) that $|\eta|^2+|v|^{2}$ and $r$ are both first order symbols belonging to the class $S(M,\Gamma)$. On the other hand, by using that the estimate
$$|r| \lesssim |\eta|^2+|v|^{2},$$
holds uniformly \wrt the parameter $\xi \in \R^n$, since $r \in S\big(|\eta|^2+|v|^{2},\Gamma\big)$; we deduce from the G\aa rding inequality (Theorem~2.5.4 in~\cite{Le}) that
$$|(r^w u, u)|  \lesssim \|D_vu\|_{L^2}^2+\|vu\|_{L^2}^2 +\|u\|_{L^2}^2.$$
Setting
\begin{equation}\label{fili5}
\Psi=\frac{|\xi|^2}{2\lambda^{4/3}}\psi\left(\frac{|\eta|^2+|v|^2}{\lambda^{2/3}}\right),
\end{equation}
we can therefore find a positive constant $C>0$ such that for all $u \in \mathcal{S}(\rr_v^n)$ and $\xi \in \rr^n$,
\begin{equation} \label{333}
-\textrm{Re}(iv.\xi u,Gu) \geq (\Psi^w u,u)-C\|D_vu\|_{L^2}^2-C\|vu\|_{L^2}^2-C\|u\|_{L^2}^2.
\end{equation}
We then deduce from (\ref{eq7fp}), (\ref{fili3}), (\ref{fili4}) and (\ref{333}) that there exists a constant
$$0 <\eps_0 \leq 1,$$
and a new positive constant $C>0$ such that for all $u \in \mathcal{S}(\rr_v^n)$ and $\xi \in \rr^n$,
\begin{equation}\label{eq7fp2}
\textrm{Re}(\langle\xi\rangle^{-s}Pu,\langle  \xi \rangle^{s}(1-\eps G)u) \geq \frac{1}{2}(\|D_vu\|_{L^2}^2+\|v u\|_{L^2}^2)+\eps_0(\Psi^wu,u)-C\|u\|_{L^2}^2.
\end{equation}
By considering separately the two regions of the phase space where, $$|\eta|^2+|v|^2 \lesssim \lambda^{2/3},$$
and $|\eta|^2+|v|^2 \gtrsim \lambda^{2/3}$; according to the support of the function
$$\psi\left(\frac{|\eta|^2+|v|^2}{\lambda^{2/3}}\right);$$
we notice that one can find a positive constant $\eps_1 >0$ such that for all $(v,\eta,\xi) \in \rr^{3n}$,
\begin{multline} \label{crucial}
\eps_0 \frac{|\xi|^2}{2\lambda^{4/3}}
\psi\left(\frac{|\eta|^2+|v|^2}{\lambda^{2/3}}\right) + \frac{1}{2}(|v|^2+|\eta|^2)
\geq \eps_1 \lambda^{2/3}+ \frac{1}{4}(|v|^2+|\eta|^2) \\ \geq \eps_1( |\xi|^{2/3} + |v|^2 +  |\eta|^2 ).
\end{multline}
This estimate is the crucial step where we combine the ellipticity in the variables $(v,\eta)$ of the real part of the symbol $p$; together with the ellipticity in the variable $\xi$ of the iterated commutator
$$[(\textrm{Im }p)^w,[(\textrm{Re }p)^w,(\textrm{Im }p)^w]]=2|\xi|^2,$$
in order to derive the optimal hypoelliptic estimate with loss of $4/3$ derivatives.
Notice from Lemma~\ref{2.2fp} and (\ref{ya2}) that
$$\eps_0 \frac{|\xi|^2}{2\lambda^{4/3}}
\psi\left(\frac{|\eta|^2+|v|^2}{\lambda^{2/3}}\right) + \frac{1}{2}(|v|^2+|\eta|^2)$$
and
$$\eps_1( |\xi|^{2/3} + |v|^2 +  |\eta|^2 ),$$
are both first order symbols belonging to the class $S(M,\Gamma)$. Recalling (\ref{fili5}) and (\ref{eq7fp2}), we can then deduce from (\ref{crucial}) and another use of the G\aa rding inequality that
there exists a new positive constant $C>0$ such that for all $s \in \rr$, $\xi \in \rr^n$ and $u \in \mathcal{S}(\rr_v^n)$,
$$\textrm{Re}(\langle\xi\rangle^{-s}Pu,\langle\xi\rangle^{s}(1-\eps G)u\big) \geq \eps_1 (\|D_vu\|_{L^2}^2+\|v u\|_{L^2}^2+\||\xi|^{1/3} u\|_{L^2}^2) - C\|u\|_{L^2}^2.$$
Notice that
$$\langle\xi\rangle^{s}(1-\eps G)=(1-\eps G)\langle\xi\rangle^{s}.$$
Recalling that the multiplier $G$ defines a bounded operator on $L^2$, Proposition~\ref{KFP2} then follows from the Cauchy-Schwarz inequality. \fin

\noindent
Taking $s=0$ in Proposition~\ref{KFP2} gives the first non optimal hypoelliptic estimate:

\bigskip

\begin{proposition}\label{KFP3}
There exists a positive constant $C > 0$ such that for all  $\xi \in \rr^n$ and $u \in \sss (\R^{n}_{v})$,
 $$\||\xi|^{1/3}u\|_{L^2}^2+\|vu\|_{L^2}^2  +  \|D_vu\|_{L^2}^2 \leq C  (\|Pu\|_{L^2}^2 + \|u\|_{L^2}^2),$$
 where $\|\cdot\|_{L^2}$ stands for the $L^2(\rr_{v}^{n})$-norm.
\end{proposition}

\bigskip

\noindent
In order to get the optimal hypoelliptic estimate, we then use an argument of commutation.

\bigskip

\begin{proposition}\label{KFP4}
There exists a positive constant $C > 0$ such that for all  $\xi \in \rr^n$ and $u \in \sss (\R^{n}_{v})$,
\begin{equation} \label{fpp}
\|\langle \xi \rangle^{2/3}u\|_{L^2}^2 + \|\langle v \rangle^2u\|_{L^2}^2 + \|\langle D_v \rangle^2u\|_{L^2}^2 \leq C(\|Pu\|_{L^2}^2 + \|u\|_{L^2}^2),
\end{equation}
where $\|\cdot\|_{L^2}$ stands for the $L^2(\rr_{v}^{n})$-norm.
\end{proposition}

\bigskip

\proof We shall successively estimate from above the three terms appearing in the left-hand-side of (\ref{fpp}).
Regarding the first one, we use Proposition~\ref{KFP2} with $s = 1/3$, to obtain that there exists a positive constant $C>0$ such that for all $\xi \in \rr^n$ and $u \in \sss (\R^{n}_{v})$,
$$\|\langle\xi\rangle^{1/3}u\|_{L^2}^2 \leq C(\|\langle\xi\rangle^{-1/3}Pu\|_{L^2}\|\langle\xi\rangle^{1/3}u\|_{L^2} +\|u\|_{L^2}^2).$$
Substituting $\langle \xi \rangle^{1/3}u$ to $u$ gives
\begin{equation}\label{fili6}
\|\langle\xi\rangle^{2/3} u\|_{L^2}^2 \leq C(\|\langle \xi \rangle^{-1/3}P\langle \xi \rangle^{1/3}u\|_{L^2 }\|\langle \xi \rangle^{2/3}u\|_{L^2} +\|\langle \xi \rangle^{1/3}u\|_{L^2}^2).
\end{equation}
Notice that
$$P=\seq{\xi}^{-1/3}P  \seq{\xi}^{1/3}.$$
It easily follows from (\ref{fili6}) that there exists a new positive constant $C>0$ such that for all $\xi \in \rr^n$ and $u \in \sss (\R^{n}_{v})$,
\begin{equation} \label{fppxi}
\|\langle \xi \rangle^{2/3}u\|_{L^2}^2  \leq C(\|Pu\|_{L^2}^2+\|u\|_{L^2}^2).
\end{equation}
For the second term, we may write
$$\|\langle v \rangle u \|_{L^2}^2=(\langle v\rangle^2 u,u) \leq \Re(Pu,u)+\|u\|_{L^2}^2 \leq \|\langle v\rangle^{-1}Pu\|_{L^2} \|\langle v\rangle u\|_{L^2}+\|u\|_{L^2}^2.$$
Substituting $\langle v\rangle u$ to $u$ gives
\begin{multline}\label{coucoumilieu}
\|\langle v\rangle^2u\|_{L^2}^2 \leq \|\langle v \rangle^{-1}P  \langle v\rangle u\|_{L^2}\|\langle v \rangle^{2}u\|_{L^2}+\|\langle v\rangle u\|_{L^2}^2 \\
 \leq \|Pu\|_{L^2}\|\langle v \rangle^{2}u\|_{L^2}+\|\langle v \rangle^{-1}[P,\langle v\rangle] u\|_{L^2}\|\langle v \rangle^{2}u\|_{L^2}+\|\langle v\rangle u\|_{L^2}^2.
\end{multline}
Symbolic calculus (Theorem~18.5.4 in~\cite{Hor85}) shows that there exist $C_b^{\infty}(\rr_v^n)$ functions $a_j$ and $b$ such that
$$[P,\langle v\rangle]=[|D_v|^2,\langle v\rangle]=\frac{1}{i}\{|\eta|^2,\langle v \rangle\}^w=\sum_{j=1}^na_j(v)D_{v_j}+b(v).$$
Here, the space $C_b^{\infty}(\rr_v^n)$ stands for the space of $C^{\infty}(\rr_v^n)$ functions whose derivatives of any order are bounded over $\rr_v^n$. It follows from Proposition~\ref{KFP3} that
\begin{equation}\label{fili7}
\|\langle v \rangle^{-1}[P,\langle v\rangle] u\|_{L^2} \leq \|[P,\langle v\rangle] u\|_{L^2} \lesssim \|D_vu\|_{L^2}+\|u\|_{L^2} \lesssim \|Pu\|_{L^2}+\|u\|_{L^2}.
\end{equation}
Finally, we easily deduce from Proposition~\ref{KFP3}, (\ref{coucoumilieu}) and (\ref{fili7}) that there exists a new positive constant $C>0$ such that for all $\xi \in \rr^n$ and $u \in \sss (\R^{n}_{v})$,
\begin{equation}\label{fili8}
\|\langle v\rangle^{2} u\|_{L^2}^2 \leq C(\|Pu\|_{L^2 }^2+\|u\|_{L^2}^2).
\end{equation}
Regarding the third term, we use similar types of estimates and write
$$\|\langle D_v \rangle u\|_{L^2}^2=(\langle D_v\rangle^2 u,u) \leq \Re(Pu,u)+\|u\|_{L^2}^2 \leq \|\langle D_v\rangle^{-1}Pu\|_{L^2} \|\langle D_v\rangle u\|_{L^2}+\|u\|_{L^2}^2.$$
Let $w$ be a $C^{\infty}(\rr,[0,1])$ function such that $w=1$ on $\rr \setminus [-2,2]$ and $w=0$ on $[-1,1]$.
Substituting $w(\langle \xi \rangle^{-1/3} \langle D_v \rangle)\langle D_v\rangle u$ to $u$ gives
\begin{multline*}
\|w(\langle \xi \rangle^{-1/3} \langle D_v \rangle)\langle D_v\rangle^2u\|_{L^2}^2 \leq  \|\langle D_v \rangle^{-1}P  w(\langle \xi \rangle^{-1/3} \langle D_v \rangle)\langle D_v\rangle u\|_{L^2}\\ \times \|w(\langle \xi \rangle^{-1/3} \langle D_v \rangle)\langle D_v\rangle^2 u\|_{L^2} + \|\langle D_v\rangle u\|_{L^2}^2,
\end{multline*}
that is
\begin{multline*}
\|w(\langle \xi \rangle^{-1/3} \langle D_v \rangle)\langle D_v\rangle^2u\|_{L^2}^2
 \leq \|w(\langle \xi \rangle^{-1/3} \langle D_v \rangle)\langle D_v\rangle^2u\|_{L^2} ( \|Pu\|_{L^2} \\ +\|\langle D_v \rangle^{-1}[P,w(\langle \xi \rangle^{-1/3} \langle D_v \rangle)\langle D_v\rangle] u\|_{L^2})+\|\langle D_v\rangle u\|_{L^2}^2.
\end{multline*}
It follows that
\begin{multline}\label{coucoumilieu4}
\|w(\langle \xi \rangle^{-1/3} \langle D_v \rangle)\langle D_v\rangle^2u\|_{L^2}^2
 \leq 2\|\langle D_v \rangle^{-1}[P,w(\langle \xi \rangle^{-1/3} \langle D_v \rangle)\langle D_v\rangle] u\|_{L^2}^2\\ +2\|Pu\|_{L^2}^2 +2\|\langle D_v\rangle u\|_{L^2}^2.
\end{multline}
Symbolic calculus (Theorem~18.5.4 in~\cite{Hor85}) shows that we have the exact identity
\begin{multline*}
[P,w(\langle \xi \rangle^{-1/3} \langle D_v \rangle)\langle D_v\rangle]=[iv.\xi+|v|^2,w(\langle \xi \rangle^{-1/3} \langle D_v \rangle)\langle D_v\rangle]\\
=\frac{1}{i}\{iv.\xi+|v|^2,w(\langle \xi \rangle^{-1/3} \langle \eta \rangle)\langle \eta \rangle\}^w.
\end{multline*}
Notice that
$$\frac{1}{i}\{iv.\xi+|v|^2,w(\langle \xi \rangle^{-1/3} \langle \eta \rangle)\langle \eta \rangle\}=(2iv-\xi).A_{\xi}(\eta),$$
with
\begin{multline*}
A_{\xi}(\eta)=\nabla_{\eta}\big(w(\langle \xi \rangle^{-1/3} \langle \eta \rangle)\langle \eta \rangle\big)=w(\langle \xi \rangle^{-1/3} \langle \eta \rangle)\nabla_{\eta}\big(\langle \eta \rangle\big)\\
+w'(\langle \xi \rangle^{-1/3} \langle \eta \rangle)\langle \xi \rangle^{-1/3} \langle \eta \rangle \nabla_{\eta}\big(\langle \eta \rangle\big),
\end{multline*}
a function satisfying
$$\forall \alpha \in \nn^n, \exists C_{\alpha}>0, \forall \xi \in \rr^n, \forall \eta \in \rr^n, \ |\partial_{\eta}^{\alpha}A_{\xi}(\eta)| \leq C_{\alpha}.$$
Symbolic calculus (Theorem~18.5.4 in~\cite{Hor85}) shows that there exists a $C_b^{\infty}(\rr_{v,\eta,\xi}^{3n})$ function $F$ such that
$$[P,w(\langle \xi \rangle^{-1/3} \langle D_v \rangle)\langle D_v\rangle]=A_{\xi}(D_v).(2iv-\xi)+F(v,\eta,\xi)^w.$$
It follows from the Calder\'on-Vaillancourt Theorem that
\begin{align*}
& \ \|\langle D_v \rangle^{-1}[P,w(\langle \xi \rangle^{-1/3} \langle D_v \rangle)\langle D_v\rangle] u\|_{L^2} \lesssim \|\langle D_v \rangle^{-1}A_{\xi}(D_v).(2iv-\xi)u\|_{L^2}+\|u\|_{L^2} \\
\lesssim & \ \|\langle v \rangle u\|_{L^2}+\|w(\langle \xi \rangle^{-1/3} \langle D_v \rangle)\langle D_v \rangle^{-1}\langle \xi \rangle u\|_{L^2}
+\|w'(\langle \xi \rangle^{-1/3} \langle D_v \rangle)\langle D_v \rangle^{-1}\langle \xi \rangle u\|_{L^2}\\
\lesssim & \  \|\langle v \rangle u\|_{L^2}+ \|\langle \xi \rangle^{2/3} u\|_{L^2},
\end{align*}
since
$$\langle \xi \rangle^{-1/3} \langle \eta \rangle \sim 1,$$
on the support of the function $w'(\langle \xi \rangle^{-1/3} \langle \eta \rangle);$ and
$$\langle \eta \rangle^{-1}\langle \xi \rangle \leq \langle \xi \rangle^{2/3},$$
on the support of the two functions
$$w(\langle \xi \rangle^{-1/3} \langle \eta \rangle) \textrm{ and } w'(\langle \xi \rangle^{-1/3} \langle \eta \rangle).$$
It follows from Proposition~\ref{KFP3}, (\ref{fppxi}) and (\ref{coucoumilieu4}) that
\begin{multline}\label{fili9}
\|w(\langle \xi \rangle^{-1/3} \langle D_v \rangle)\langle D_v\rangle^2u\|_{L^2}^2 \lesssim \|Pu\|_{L^2}^2 +\|\langle D_v\rangle u\|_{L^2}^2\\ +\|\langle v \rangle u\|_{L^2}^2+ \|\langle \xi \rangle^{2/3} u\|_{L^2}^2
\lesssim  \|Pu\|_{L^2}^2 +\|u\|_{L^2}^2.
\end{multline}
Notice that there exists a positive constant $C_0>0$ such that for all $\xi \in \rr^n$ and $\eta \in \rr^n$,
\begin{equation}\label{fili10}
\frac{1}{C_0}\big(\langle \xi \rangle^{4/3}+\langle \eta \rangle^{4}\big) \leq \langle \xi \rangle^{4/3}+w(\langle \xi \rangle^{-1/3} \langle \eta \rangle)^2\langle \eta \rangle^{4},
\end{equation}
because $\langle \eta \rangle^4 \leq 2 \langle \xi \rangle^{4/3}$, when $w(\langle \xi \rangle^{-1/3} \langle \eta \rangle) \neq 1$.
Finally, by collecting the estimates (\ref{fppxi}), (\ref{fili8}), (\ref{fili9}) and (\ref{fili10}), we find that there exists a new positive constant $C > 0$ such that for all  $\xi \in \rr^n$ and $u \in \sss (\R^{n}_{v})$,
$$\|\langle \xi \rangle^{2/3}u\|_{L^2}^2 + \|\langle v \rangle^2u\|_{L^2}^2 + \|\langle D_v \rangle^2u\|_{L^2}^2 \leq C(\|Pu\|_{L^2}^2 + \|u\|_{L^2}^2),$$
which proves Proposition~\ref{KFP4}.
\fin

\bigskip

\noindent
When coming back to the direct side and integrating with respect to the $x$ variable, Proposition~\ref{KFP} directly follows from Proposition \ref{KFP4}. This proves the optimal hypoelliptic estimate fulfilled by the Fokker-Planck operator without external potential. \fin

%%%%%%%%%%%%%%%%%%%%%%%%%%%%%%%%%%%%%%%%%%%%%%%%%%%%%%%%%%%%%
\section{Anisotropic hypoelliptic estimates for linear Landau-type operators}\label{la1}

In this section, we consider the class of linear Landau-type operators
\begin{equation}\label{ind1Q}
P= iv.D_x+ D_v.\lambda(v)D_{v}+(v \wedge D_v).\mu(v) (v \wedge D_v)+F(v), \ x,v \in \rr^{3};
\end{equation}
that is
\begin{multline*}
P = i\sum_{j=1}^3v_jD_{x_j}+ \sum_{j=1}^3D_{v_j}\lambda(v)D_{v_j}+(v_2D_{v_3}-v_3D_{v_2})\mu(v)(v_2D_{v_3}-v_3D_{v_2})\\
+(v_3D_{v_1}-v_1D_{v_3})\mu(v)(v_3D_{v_1}-v_1D_{v_3})+(v_1D_{v_2}-v_2D_{v_1})\mu(v)(v_1D_{v_2}-v_2D_{v_1})+F(v),
\end{multline*}
with $D_x=i^{-1}\partial_x$, $D_v=i^{-1}\partial_v$ and $\gamma \in [-3,1]$; where the diffusion is given by smooth positive functions $\lambda$, $\mu$ and $F$ satisfying for all $\alpha \in \nn^3$, there exists $C_{\alpha}>0$ such that
\begin{equation}\label{eq-3.1Q}
 \forall v \in \rr^3,  |\partial_v^{\alpha}\lambda(v)| +|\partial_v^{\alpha}\mu(v)| \leq C_{\alpha} \langle v \rangle^{\gamma-|\alpha|}; |\partial_v^{\alpha}F(v)| \leq C_{\alpha} \langle v \rangle^{\gamma+2-|\alpha|};
\end{equation}
and
\begin{equation}\label{w2Q}
\exists C>0, \forall v \in \rr^3,  \ \lambda(v) \geq C \langle v \rangle^{\gamma}; \ \mu(v) \geq C \langle v \rangle^{\gamma}; \  F(v) \geq C \langle v \rangle^{\gamma+2};
\end{equation}
with $\langle v \rangle=(1+|v|^2)^{1/2}$. We aim at proving the optimal anisotropic hypoelliptic estimate with loss of $4/3$ derivatives given in Theorem~\ref{mainlandau}.

In order to do so, we begin by considering generalized linear Landau-type operators
\begin{equation}\label{fil1}
P=iv.D_x+ \sum_{j,k=1}^nD_{v_j}A_{j,k}(v)D_{v_k}+F(v);
\end{equation}
where $x,v \in \rr^{n}$,
 $D_x=i^{-1}\partial_x$, $D_v=i^{-1}\partial_v$, $\gamma \in [-3,1]$.
  Here $A(v)=(A_{j,k}(v))_{1 \leq j,k \leq n}$ stands for a positive definite symmetric matrix with real-valued smooth entries verifying
\begin{equation}\label{eq-3QQ}
|\partial_v^{\alpha}A_{j,k}(v)| \lesssim \langle v \rangle^{\gamma+2-|\alpha|}, \ \alpha \in \nn^{n}, \ 1 \leq j,k \leq n;
\end{equation}
and $F$ is a smooth positive function verifying
\begin{equation}\label{eq-3.1QQ}
 F(v) \gtrsim \langle v \rangle^{\gamma+2} \textrm{ and } |\partial_v^{\alpha}F(v)| \lesssim \langle v \rangle^{\gamma+2-|\alpha|}, \ \alpha \in \nn^{n}.
\end{equation}
We recall that the notation
$$f(v) \lesssim g(v),$$
means that there exists a positive constant $C>0$ such that the estimate
$$f(v) \leq C g(v),$$
is fulfilled for all $v \in \rr^{n}$.
We assume that we may write
\begin{equation}\label{eq-2QQ}
A(v)=B(v)^TB(v),
\end{equation}
where $B(v)$ is a  matrix with real-valued smooth entries verifying
\begin{equation}\label{eq-1QQ}
|\partial_v^{\alpha}B_{j,k}(v)| \lesssim \langle v \rangle^{\frac{\gamma}{2}+1-|\alpha|},\ \alpha \in \nn^{n}, \ 1 \leq j,k \leq n;
\end{equation}
and $B(v)^T$ is its adjoint.
Moreover, we assume that there exists a constant $c>0$ such that for all $v, \eta \in \rr^n$,
\begin{equation}\label{eq0QQ}
A(v)\eta.\eta =|B(v)\eta|^2\geq c\langle v \rangle^{\gamma}|\eta|^2.
\end{equation}
Notice that linear Landau-type operators are particular generalized linear Landau-type operators when taking
\begin{equation}\label{an2QQ}
B(v)=\left(
\begin{array}{ccc}
\sqrt{\lambda(v)} & -v_3\sqrt{\mu(v)}& v_2\sqrt{\mu(v)}\\
v_3\sqrt{\mu(v)} &  \sqrt{\lambda(v)} & -v_1\sqrt{\mu(v)}\\
-v_2\sqrt{\mu(v)}& v_1\sqrt{\mu(v)}& \sqrt{\lambda(v)}
\end{array}
\right),
\end{equation}
with $\lambda $ and $\mu$ being the functions defined in (\ref{eq-3.1Q}) and (\ref{w2Q}).
Indeed, we have for any $\eta \in \rr^3$,
\begin{equation}\label{rih1}
|B(v)\eta|^2=|\sqrt{\lambda(v)}\eta+\sqrt{\mu(v)}v \wedge \eta|^2=|\sqrt{\lambda(v)}\eta|^2+
|\sqrt{\mu(v)}v \wedge \eta|^2
\geq c \langle v \rangle^{\gamma}|\eta|^2.
\end{equation}

\subsection{First estimates for generalized linear Landau-type operators}

In order to prove Theorem~\ref{mainlandau}, we shall use a multiplier method inspired from the one presented in the previous section for the Fokker-Planck operator without external potential. Recalling (\ref{eq-2QQ}), the Weyl symbol of a generalized linear landau-type operator (\ref{fil1}) may write as
$$iv.\xi+|B(v)\eta|^2+F(v)+\textrm{\textsl{Lower order terms}}.$$
By denoting
$$\tilde{p}=iv.\xi+|B(v)\eta|^2+F(v),$$
we shall take advantage of the ellipticity in the variables $(v,\eta)$ of the real part of the symbol $\tilde{p}$,
$$\textrm{Re }\tilde{p}=|B(v)\eta|^2+F(v).$$
As in the case of the Fokker-Planck operator, the main point in proving Theorem~\ref{mainlandau} is then to get a control of the $\xi$ variable. Notice again that this control cannot be derived from the ellipticity of the symbol $\tilde{p}$; and that we will need to consider the following iterated commutator
$$[(\textrm{Im }\tilde{p})^w,[(\textrm{Re }\tilde{p})^w,(\textrm{Im }\tilde{p})^w]],$$
where $\textrm{Re }\tilde{p}$ and $\textrm{Im }\tilde{p}$ stand for the real and imaginary parts of the symbol $\tilde{p}$; in order to get some ellipticity in the $\xi$ variable. Indeed, usual symbolic calculus (see Theorem~18.5.4 in~\cite{Hor85}) or a direct computation shows that the Weyl symbol of this iterated commutator is exactly given by the iterated Poisson brackets
$$-\{\textrm{Im }\tilde{p},\{\textrm{Re }\tilde{p},\textrm{Im }\tilde{p}\}\}=\{\textrm{Im }\tilde{p},\{\textrm{Im }\tilde{p},\textrm{Re }\tilde{p}\}\}=2|B(v)\xi|^2.$$
The structure of this iterated poisson bracket suggests to introduce the following anisotropic symbol
\begin{equation}\label{m2}
\lambda=\big(1+|B(v)\xi|^2+|B(v)\eta|^2+F(v)\big)^{1/2},
\end{equation}
which defines an anisotropic Sobolev scale which is exactly related to the anisotropy of the diffusion. As in the case of the Fokker-Planck operator, we aim at establishing an optimal hypoelliptic estimate with loss of $4/3$ derivatives in this anisotropic Sobolev scale
$$\|(\lambda^{2/3})^w u\|_{L^2}^2 \lesssim \|Pu\|_{L^2}^2+\|u\|_{L^2}^2.$$
By noticing that for a generalized linear Landau-type operator
$$H_{\textrm{Re}\tilde{p}}\ \textrm{Im }\tilde{p}=\{\textrm{Re }\tilde{p},\textrm{Im }\tilde{p}\}=2B(v)\xi.B(v)\eta,$$
it is natural to consider the following multiplier:
Let $\Psi$ be a $C_0^{\infty}(\rr,[0,1])$ function such that
\begin{equation}\label{n5}
\psi=1 \textrm{ on } [-1,1], \textrm{ and } \textrm{supp }\psi \subset [-2,2].
\end{equation}
Define the real-valued symbol
\begin{equation}\label{m1}
g=-\frac{B(v)\xi.B(v)\eta}{\lambda^{4/3}}\psi\left(\frac{|B(v)\eta|^2+F(v)}{\lambda^{2/3}}\right),
\end{equation}
where $\lambda$ is the symbol defined in (\ref{m2}). The main difference with the Fokker-Planck case is that this multiplier does not belong anymore to a symbol class with good symbolic calculus. Indeed, because of the anisotropy of the symbol $\tilde{p}$, we will have to deal with gainless symbolic calculus. As a consequence, the implementation of the method developed for the Fokker-Planck operator will be more complex and will require more advanced microlocal analysis. In order to handle this setting with gainless symbolic calculus, we shall use some elements of Wick calculus developed by N.~Lerner in~\cite{Ler03}. For convenience of reading, the main features and the definition of Wick calculus is recalled in a short self-contained presentation given in appendix (Section~\ref{la3}).

When studying generalized linear Landau-type operators, it is convenient to perform a partial Fourier transform in the $x$ variable; and to study these operators on the Fourier side
$$P=iv.\xi+ \sum_{j,k=1}^nD_{v_j}A_{j,k}(v)D_{v_k}+F(v), \ v,\xi \in \rr^{n};$$
where the variable $\xi$ can now be seen as a parameter. In the following, we shall therefore consider quantizations of symbols only in the variable $v$ and its dual variable $\eta$; and denote by $\|\cdot\|_{L^2}$ the $L^2(\rr_{v}^{n})$-norm. A key step in proving Theorem~\ref{mainlandau} is the proof of the following proposition somehow equivalent to Proposition~\ref{KFP2} in the Fokker-Planck case.

\bigskip

\begin{proposition}\label{prop1}
Let $s \in \rr$ and $P$ be a generalized linear Landau-type operator fulfilling the assumptions \emph{(\ref{eq-3QQ})}, \emph{(\ref{eq-3.1QQ})}, \emph{(\ref{eq-2QQ})}, \emph{(\ref{eq-1QQ})} and \emph{(\ref{eq0QQ})}. Then, there exists a constant $C>0$ such that for all $\xi \in \rr^n$ and $u \in \mathcal{S}(\rr_{v}^{n})$,
\begin{multline*}
\|m(v,\xi)u\|_{L^2}^2 +\|B(v)\nabla_vu\|_{L^2}^2+\|\sqrt{F(v)}u\|_{L^2}^2 \\ \leq C\|\langle B(v) \xi \rangle^{-s}Pu\|_{L^2}
\|\langle B(v) \xi \rangle^{s}u\|_{L^2}+C\|u\|_{L^2}^2,
\end{multline*}
where $\|\cdot\|_{L^2}$ stands for the $L^2(\rr_{v}^{n})$-norm and
$$m(v,\xi)=\Big(\int_{\rr^n}\langle B(v+\tilde{v})\xi \rangle^{2/3}e^{-2\pi\tilde{v}^2}2^{n/2}d\tilde{v}\Big)^{1/2}.$$
\end{proposition}

\bigskip

\noindent
\textit{Remark.} The use of Wick calculus accounts for the definition of the quantity $m(v,\xi)$ in Proposition~\ref{prop1}. We shall see that the function $m(v,\xi)^2$ is actually the Wick quantization of the symbol $\langle B(v)\xi \rangle^{2/3}$.

\bigskip

Let $m\geq 1$  be a $C^{\infty}$ order function on $\real^{2n}$, we denote by $S(m,dv^2+d\eta^2)$ the symbol class
$$\left\{ a\in C^{\infty}(\real^{2n},\cc): \forall \alpha \in \nn^{2n}, \exists C_{\alpha}>0, \forall (v,\eta) \in \rr^{2n}, \  |\partial_{v,\eta}^{\alpha} a(v,\eta)| \leq C_{\alpha} m(v,\eta)\right\}.$$
As a starting point in the proof of Proposition~\ref{prop1}, we notice that a generalized linear Landau-type operator is accretive
$$\textrm{Re}\big(\langle B(v) \xi \rangle^{-s}Pu,\langle B(v) \xi \rangle^{s}u\big)=\textrm{Re}(Pu,u)=\|B(v)\nabla_vu\|_{L^2}^2+\big\|\sqrt{F(v)}u\big\|_{L^2}^2 \geq 0,$$
for any $s \in \rr$.
It follows from the Cauchy-Schwarz inequality and (\ref{eq0QQ}) that
\begin{multline}\label{eq1}
\|\langle v \rangle^{\gamma/2}\nabla_vu\|_{L^2}^2+\big\|\sqrt{F(v)}u\big\|_{L^2}^2 \lesssim \|B(v)\nabla_vu\|_{L^2}^2+\big\|\sqrt{F(v)}u\big\|_{L^2}^2  \\ \leq \|\langle B(v) \xi \rangle^{-s}Pu\|_{L^2}
\|\langle B(v) \xi \rangle^{s}u\|_{L^2},
\end{multline}
uniformly with respect to the parameter $\xi \in \rr^n$.

\bigskip

\begin{lemma}\label{lem0}
For $m \in \rr$,
$$\langle B(v) \xi \rangle^m \in S\big(\langle B(v) \xi \rangle^m,dv^2+d\eta^2\big) \textrm{ and } \lambda^m \in S(\lambda^m,dv^2+d\eta^2),$$
uniformly with respect to the parameter $\xi$ in $\rr^n$.

\end{lemma}

\bigskip

\noindent
\textit{Proof of Lemma~\ref{lem0}.} Notice from (\ref{eq-3.1QQ}), (\ref{eq-1QQ}) and (\ref{eq0QQ}) that
$$|\partial_v^{\alpha}F(v)| \lesssim F(v), \ \alpha \in \nn^n,$$
$$|\partial_v^{\alpha}B(v)|^2 \lesssim F(v), \ \alpha \in \nn^n,$$
\begin{equation}\label{n1b}
|\partial_v^{\alpha}B(v)\xi|^2 \lesssim \langle v \rangle^{\gamma}|\xi|^2 \lesssim
|B(v)\xi|^2,
\end{equation}
when $\alpha \in \nn^n$ with $|\alpha| \geq 1$; and
\begin{equation}\label{n1}
|\partial_v^{\beta}B(v)\eta|^2 \lesssim \langle v \rangle^{\gamma}|\eta|^2 \lesssim |B(v)\eta|^2,
\end{equation}
when $\beta \in \nn^n$ with $|\beta| \geq 1$. One can then deduce by using the Cauchy-Schwarz inequality and these estimates that for all $\alpha, \beta \in \nn^n$,
$$\big|\partial_{v}^{\alpha}\big(\langle B(v) \xi \rangle^2\big)\big| \lesssim \langle B(v) \xi \rangle^2 \textrm{ and } |\partial_{v}^{\alpha}\partial_{\eta}^{\beta}(\lambda^2)| \lesssim \lambda^2;$$
uniformly with respect to the parameter $\xi$ in $\rr^n$. Lemma~\ref{lem0} directly follows from those estimates. $\Box$

\bigskip

\begin{lemma}\label{lem-1}
We have
$$\psi\left(\frac{|B(v)\eta|^2+F(v)}{\lambda^{2/3}}\right) \in S(1,dv^2+d\eta^2),$$
uniformly with respect to the parameter $\xi$ in $\rr^n$.

\end{lemma}

\bigskip

\noindent
\textit{Proof of Lemma~\ref{lem-1}.} Notice from (\ref{eq-3.1QQ}), (\ref{eq-1QQ}), (\ref{eq0QQ}), (\ref{n5}) and (\ref{n1}) that on the support of the function
$$\psi\left(\frac{|B(v)\eta|^2+F(v)}{\lambda^{2/3}}\right),$$
we have
\begin{multline}\label{n2}
|\partial_v^{\alpha_1}B(v)\eta|^2 +|\partial_v^{\alpha_2}B(v)|^2+|\partial_v^{\alpha_3}F(v)| \lesssim \langle v \rangle^{\gamma}|\eta|^2 + \langle v \rangle^{\gamma+2}\\ \lesssim |B(v)\eta|^2+ F(v) \lesssim \lambda^{2/3},
\end{multline}
when $\alpha_1, \alpha_2, \alpha_3 \in \nn^n$ with $|\alpha_1| \geq 1$. It follows from the Cauchy-Schwarz inequality that any derivatives of the term $|B(v)\eta|^2+F(v)$ can be estimated from above by a constant times the term $\lambda^{2/3}$ on the support of this function. One can therefore directly deduce the result of Lemma~\ref{lem-1} from Lemma~\ref{lem0}. $\Box$

\bigskip

\begin{lemma}\label{lem1}
The symbol $g$ belongs to the class $S(1,dv^2+d\eta^2)$ uniformly with respect to the parameter $\xi$ in $\rr^n$.
\end{lemma}

\bigskip

\noindent
\textit{Proof of Lemma~\ref{lem1}.} Notice from (\ref{n5}) that
\begin{equation}\label{m3}
|B(v)\eta|^2+F(v) \leq 2\lambda^{2/3},
\end{equation}
on the support of the function
$$\psi\left(\frac{|B(v)\eta|^2+F(v)}{\lambda^{2/3}}\right).$$
By recalling (\ref{m2}) and using that $|B(v)\xi| \leq \lambda$, we deduce from the Cauchy-Schwarz inequality that one can estimate
\begin{equation}\label{n6}
|B(v)\xi.B(v)\eta| \leq |B(v)\xi||B(v)\eta| \leq \sqrt{2}\lambda^{4/3},
\end{equation}
on this support. The symbol $g$ is therefore a bounded function uniformly with respect to the parameter $\xi$ in $\rr^n$. We saw in (\ref{n1b}) that one can always estimate from above
\begin{equation}\label{an1}
|\partial_v^{\alpha}B(v)\xi| \lesssim |B(v)\xi| \lesssim \lambda.
\end{equation}
Since from (\ref{n2}), one can estimate from above the modulus of all the derivatives of the term $B(v)\eta$ by a constant times $\lambda^{1/3}$ on the support of the function
$$\psi\left(\frac{|B(v)\eta|^2+F(v)}{\lambda^{2/3}}\right),$$
it follows from the Cauchy-Schwarz inequality and (\ref{an1}) that one can estimate from above the modulus of all the derivatives of the term $B(v)\xi.B(v)\eta$ by a constant times $\lambda^{4/3}$ on this support. According to Lemma~\ref{lem0} and Lemma~\ref{lem-1}, this proves that the symbol $g$ belongs to the class $S(1,dv^2+d\eta^2)$ uniformly with respect to the parameter $\xi$ in $\rr^n$; and ends the proof of Lemma~\ref{lem1}. $\Box$

\subsubsection{Some symbolic calculus}
We shall consider the multiplier $G=g^{\textrm{Wick}}$ defined by the Wick quantization of the symbol $g$. We refer the reader to the appendix on Wick calculus at the end of this note for the definition of this quantization and a recall of its main features.

We begin by noticing from (\ref{lay1bis}) and (\ref{lay2bis}) that there exists a real-valued symbol $\tilde{g}$ belonging to the class $S(1,dv^2+d\eta^2)$ uniformly with respect to the parameter $\xi$ in $\rr^n$ such that
\begin{equation}\label{m4}
G=g^{\textrm{Wick}}=\tilde{g}^w;
\end{equation}
where $\tilde{g}^w$ denotes the operator obtained by the Weyl quantization of the symbol $\tilde{g}$ with the normalization
\begin{equation}\label{weyl1}
\tilde{g}^wu(v)=\frac{1}{(2\pi)^n}\int_{\rr^{2n}}e^{i(v-\tilde{v}).\eta}\tilde{g}\Big(\frac{v+\tilde{v}}{2},\eta\Big)u(\tilde{v})d\tilde{v}d\eta.
\end{equation}

\bigskip

\begin{lemma}\label{lem2}
If $a \in S(1,dv^2+d\eta^2)$ then there exists $c_1>0$ such that for all $u \in \mathcal{S}(\rr^n)$,
$$\big\|\big[a^w,\sqrt{F(v)}\big]u\big\|_{L^2} \leq c_1 \big\|\sqrt{F(v)}u\big\|_{L^2},$$
where $\big[a^w,\sqrt{F(v)}\big]$ denotes the commutator of the operators $a^w$ and $\sqrt{F(v)}$.
\end{lemma}

\bigskip

\noindent
\textit{Proof of Lemma~\ref{lem2}.} Notice from (\ref{eq-3.1QQ}) that
$$F(v)^m \in S\big(F(v)^m,dv^2+d\eta^2\big),$$
for any $m \in \rr$; and that from symbolic calculus the Weyl symbol of the operator
$$\big[a^w,\sqrt{F(v)}\big]\big(\sqrt{F(v)}\big)^{-1},$$
therefore belongs to the symbol class $S(1,dv^2+d\eta^2)$. Lemma~\ref{lem2} then directly follows from the Calder\'on-Vaillancourt Theorem
$$\big\|\big[a^w,\sqrt{F(v)}\big]u\big\|_{L^2} =\big\|\big[a^w,\sqrt{F(v)}\big]\big(\sqrt{F(v)}\big)^{-1}\sqrt{F(v)} u\big\|_{L^2} \lesssim \big\|\sqrt{F(v)}u\big\|_{L^2}. \quad \Box$$

\bigskip

\begin{lemma}\label{lem3}
We have
$$\big\|B(v)\nabla_v(Gu)\big\|_{L^2}^2 \lesssim  \|B(v)\nabla_vu\|_{L^2}^2+ \big\|\sqrt{F(v)}u\big\|_{L^2}^2,$$
uniformly with respect to the parameter $\xi$ in $\rr^n$.
\end{lemma}

\bigskip

\noindent
\textit{Proof of Lemma~\ref{lem3}.}
Recalling from (\ref{m4}) that $G=\mathcal{O}_{\mathcal{L}(L^2)}(1)$, since
$$\tilde{g} \in S(1,dv^2+d\eta^2),$$
together with (\ref{eq1}); we notice that
it is sufficient to prove that
\begin{equation}\label{eq3.5}
\big\|\big[b(v)\nabla_v,G\big]u\big\|_{L^2}^2 \lesssim \|\langle v \rangle^{\gamma/2}\nabla_vu\|_{L^2}^2+ \big\|\sqrt{F(v)}u\big\|_{L^2}^2,
\end{equation}
when $b$ is a smooth function fulfilling the estimates (\ref{eq-1QQ}).
By writing that
\begin{equation}\label{eq4}
\big\|\big[b(v)\nabla_v,G\big]u\big\|_{L^2}^2  \leq 2\big\|b(v)\big[\nabla_v,G\big]u\big\|_{L^2}^2 +2\big\|\big[b(v),G\big] \nabla_v u\big\|_{L^2}^2,
\end{equation}
we notice from (\ref{eq-3.1QQ}), (\ref{eq-1QQ}), (\ref{m4}), the Calder\'on-Vaillancourt Theorem and Lemma~\ref{lem2} that
\begin{align}
& \ \big\|b(v)\big[\nabla_v,G\big]u\big\|_{L^2} \lesssim   \big\|b(v)(\nabla_{v}\tilde{g})^wu\big\|_{L^2} \\
\lesssim & \ \big\| \langle v \rangle^{\frac{\gamma}{2}+1}(\nabla_{v}\tilde{g})^wu\big\|_{L^2}  \lesssim  \big\| \sqrt{F(v)}(\nabla_{v}\tilde{g})^wu\big\|_{L^2}\\
 \lesssim & \ \big\| \big[\sqrt{F(v)},(\nabla_{v}\tilde{g})^w\big]u\big\|_{L^2}+\big\|(\nabla_{v}\tilde{g})^w \sqrt{F(v)}u\big\|_{L^2} \lesssim \big\|\sqrt{F(v)}u\big\|_{L^2}.\label{eq5}
\end{align}
Recalling that $\tilde{g} \in S(1,dv^2+d\eta^2)$ together with (\ref{eq-1QQ}) and (\ref{m4}), symbolic calculus (Theorem~2.3.8 and Corollary~2.3.10 in~\cite{Le}) ensures that the Weyl symbol of the operator
$$\big[b(v),G\big]\langle v \rangle^{-\gamma/2},$$
belongs to the class $S(1,dv^2+d\eta^2)$, and it follows from the Calder\'on-Vaillancourt Theorem that
$$\big\|\big[b(v),G\big] \nabla_v u\big\|_{L^2} =\big\|\big[b(v),G\big]\langle v \rangle^{-\gamma/2}\langle v \rangle^{\gamma/2} \nabla_v u\big\|_{L^2} \lesssim \big\|\langle v \rangle^{\gamma/2} \nabla_v u\big\|_{L^2},$$
which together with (\ref{eq4}) and (\ref{eq5}) proves the estimate (\ref{eq3.5}) and ends the proof of Lemma~\ref{lem3}. $\Box$

\bigskip

\begin{lemma}\label{lem4}
We have
$$\big|\big(F(v) u,Gu\big)\big|+\Big|\Big(\sum_{j,k=1}^nD_{v_j}A_{j,k}(v)D_{v_k}u,Gu\Big)\Big| \lesssim \|\langle B(v) \xi \rangle^{-s}Pu\|_{L^2}
\|\langle B(v) \xi \rangle^{s}u\|_{L^2},$$
uniformly with respect to the parameter $\xi \in \rr^n$.
\end{lemma}

\bigskip

\noindent
\textit{Proof of Lemma~\ref{lem4}.}
We may write
\begin{multline}\label{m5}
\big(F(v) u,Gu\big)=\big(GF(v) u,u\big)=
\big(G\sqrt{F(v)} u,\sqrt{F(v)}u\big)\\ +\big(\big[G,\sqrt{F(v)}\big]\sqrt{F(v)}u, u\big) =
\big(G\sqrt{F(v)} u,\sqrt{F(v)}u)+(\sqrt{F(v)}u,\big[\sqrt{F(v)},G\big] u\big),
\end{multline}
since the operator $G$ whose Weyl symbol is real-valued is formally selfadjoint on $L^2$.
Recalling from (\ref{m4}) that $G=\mathcal{O}_{\mathcal{L}(L^2)}(1)$, since $\tilde{g} \in S(1,dv^2+d\eta^2)$; we deduce from (\ref{m4}), (\ref{m5}), the triangle inequality, the Cauchy-Schwarz inequality, the Calder\'on-Vaillancourt Theorem and Lemma~\ref{lem2} that
$$\big|\big(F(v) u,Gu\big)\big| \lesssim \big\|\sqrt{F(v)}u\big\|_{L^2}^2+\big\|\big[G,\sqrt{F(v)}\big]u\big\|_{L^2}^2  \lesssim \big\|\sqrt{F(v)}u\big\|_{L^2}^2,$$
which implies by using (\ref{eq1}) that
\begin{equation}\label{eq2}
\big|\big(F(v) u,Gu\big)\big|  \lesssim \|\langle B(v) \xi \rangle^{-s}Pu\|_{L^2}
\|\langle B(v) \xi \rangle^{s}u\|_{L^2},
\end{equation}
uniformly with respect to the parameter $\xi \in \rr^n$.
Let us now notice from (\ref{eq-2QQ}), (\ref{eq1}) and Lemma~\ref{lem3} that
\begin{align*}
& \ \Big|\Big(\sum_{j,k=1}^nD_{v_j}A_{j,k}(v)D_{v_k}u,Gu\Big)\Big|=\big|\big(B(v)\nabla_v u ,B(v)\nabla_v(Gu)\big)\big| \\
\leq & \ \|B(v)\nabla_v u\|_{L^2}\| B(v)\nabla_v(Gu)\|_{L^2}
\lesssim  \|B(v)\nabla_vu\|_{L^2}^2+ \big\|\sqrt{F(v)}u\big\|_{L^2}^2 \\
\lesssim & \  \|\langle B(v) \xi \rangle^{-s}Pu\|_{L^2}
\|\langle B(v) \xi \rangle^{s}u\|_{L^2},
\end{align*}
uniformly with respect to the parameter $\xi \in \rr^n$; which together with the estimate (\ref{eq2}) proves Lemma~\ref{lem4}. $\Box$

\bigskip

\noindent
Let $\eps$ be a positive parameter such that $0<\eps \leq 1$.  We use a multiplier method and write that
\begin{multline}\label{eq7}
\textrm{Re}\big(\langle B(v) \xi \rangle^{-s}Pu,\langle B(v) \xi \rangle^{s}(1-\eps G)u\big)=\|B(v)\nabla_vu\|_{L^2}^2+\big\|\sqrt{F(v)}u\big\|_{L^2}^2\\ -\eps \textrm{Re}\big(iv.\xi u,Gu\big)
-\eps \textrm{Re}\Big(\sum_{j,k=1}^nD_{v_j}A_{j,k}(v)D_{v_k}u,Gu\Big)-\eps \textrm{Re}\big(F(v) u,Gu\big),
\end{multline}
for any $0<\eps \leq 1$.

\bigskip

\begin{lemma}\label{lem-10}
We have
$$\|\langle B(v) \xi \rangle^{s}(1-\eps G)u\|_{L^2} \lesssim \|\langle B(v) \xi \rangle^{s}u\|_{L^2},$$
uniformly with respect to the parameter $\xi \in \rr^n$.

\end{lemma}

\bigskip

\noindent
\textit{Proof of Lemma~\ref{lem-10}.}
Recalling from (\ref{m4}) that $G=\tilde{g}^w$, with $\tilde{g} \in S(1,dv^2+d\eta^2)$; we can then deduce from symbolic calculus that the Weyl symbol of the operator
$$\langle B(v) \xi \rangle^{s}(1-\eps G)\langle B(v) \xi \rangle^{-s},$$
belongs to the symbol class $S(1,dv^2+d\eta^2)$, since we know from Lemma~\ref{lem0} that
$$\langle B(v) \xi \rangle^{s} \in S\big(\langle B(v) \xi \rangle^{s},dv^2+d\eta^2\big) \textrm{ and } \langle B(v) \xi \rangle^{-s} \in S\big(\langle B(v) \xi \rangle^{-s},dv^2+d\eta^2\big).$$
Lemma~\ref{lem-10} then directly follows from the Calder\'on-Vaillancourt Theorem
\begin{align*}
\|\langle B(v) \xi \rangle^{s}(1-\eps G)u\|_{L^2}& =\|\langle B(v) \xi \rangle^{s}(1-\eps G)\langle B(v) \xi \rangle^{-s}\langle B(v) \xi \rangle^{s}u\|_{L^2} \\
 & \lesssim \|\langle B(v) \xi \rangle^{s}u\|_{L^2}.
 \end{align*}
The proof is complete. \fin

\bigskip

\noindent
One can then deduce from the Cauchy-Schwarz inequality, Lemma~\ref{lem4}, Lemma~\ref{lem-10} and (\ref{eq7}) that there exists $C>0$ such that for all $0<\eps \leq 1$ and $u \in \mathcal{S}(\rr^n)$,
\begin{multline}\label{eq8}
\|B(v)\nabla_vu\|_{L^2}^2+\big\|\sqrt{F(v)}u\big\|_{L^2}^2-\eps \textrm{Re}\big(iv.\xi u,Gu\big)  \\ \leq C\|\langle B(v) \xi \rangle^{-s}Pu\|_{L^2}
\|\langle B(v) \xi \rangle^{s}u\|_{L^2},
\end{multline}
uniformly with respect to the parameter $\xi \in \rr^n$.
Recalling (\ref{m4}) and noticing from (\ref{lay1}) and (\ref{lay2}) that $v^\textrm{Wick}=v$, we may rewrite (\ref{eq8}) as
\begin{multline}\label{eq9}
\|B(v)\nabla_vu\|_{L^2}^2+\big\|\sqrt{F(v)}u\big\|_{L^2}^2-\eps \textrm{Re}\big(i\xi.v^{\textrm{Wick}} u,g^{\textrm{Wick}}u\big) \\
 \leq C\|\langle B(v) \xi \rangle^{-s}Pu\|_{L^2}
\|\langle B(v) \xi \rangle^{s}u\|_{L^2}.
\end{multline}
By using that real Hamiltonians get quantized in the Wick quantization by formally selfadjoint operators on $L^2$, we deduce from Lemma~\ref{lem1} and (\ref{lay4}) that
\begin{align}\label{n8}
-\eps \textrm{Re}\big(i\xi.v^{\textrm{Wick}} u,g^{\textrm{Wick}}u\big)=& \ -\eps\big(\textrm{Re}\big(g^{\textrm{Wick}}(i\xi.v)^{\textrm{Wick}}\big) u,u\big)\\
=& \ \eps\frac{1}{4\pi}\big(\{\xi.v,g\}^{\textrm{Wick}}u,u\big).
\end{align}
A direct computation of the Poisson bracket using (\ref{m1}) gives that
\begin{multline}\label{n7}
\{\xi.v,g\}=\big(B(v)\xi.B(v)\eta\big)
\big(\xi.\partial_{\eta}\big(\lambda^{-4/3}\big)\big)
\psi\left(\frac{|B(v)\eta|^2+F(v)}{\lambda^{2/3}}\right)\\
+\frac{|B(v)\xi|^2}{\lambda^{4/3}}
\psi\left(\frac{|B(v)\eta|^2+F(v)}{\lambda^{2/3}}\right)
+\frac{B(v)\xi.B(v)\eta}{\lambda^{4/3}}
\xi.\partial_{\eta}\left[\psi\left(\frac{|B(v)\eta|^2+F(v)}{\lambda^{2/3}}\right)\right].
\end{multline}

\bigskip

\begin{lemma}\label{lem5}
For $m \in \rr$, we have $|\xi.\partial_{\eta}(\lambda^m)| \lesssim \lambda^m,$
uniformly with respect to the parameter $\xi$ in $\rr^n$.
\end{lemma}

\bigskip

\noindent
\textit{Proof of Lemma~\ref{lem5}.}
Lemma~\ref{lem5} follows directly from (\ref{m2}) and the fact that
$$|\xi.\partial_{\eta}(\lambda^2)|=2|B(v)\xi.B(v)\eta| \leq 2|B(v)\xi||B(v)\eta| \leq |B(v)\xi|^2+|B(v)\eta|^2 \leq \lambda^2. \quad \Box$$

\bigskip

\begin{lemma}\label{lem6}
We have
$$\left|\xi.\partial_{\eta}\left[\psi\left(\frac{|B(v)\eta|^2+F(v)}{\lambda^{2/3}}\right)\right]\right| \lesssim 1+|B(v)\eta|^2+F(v),$$
uniformly with respect to the parameter $\xi$ in $\rr^n$.
\end{lemma}

\bigskip

\noindent
\textit{Proof of Lemma~\ref{lem6}.} We may write
\begin{multline*}
\xi.\partial_{\eta}\left[\psi\left(\frac{|B(v)\eta|^2+F(v)}{\lambda^{2/3}}\right)\right]=
\psi'\left(\frac{|B(v)\eta|^2+F(v)}{\lambda^{2/3}}\right)\Big[\frac{2B(v)\eta.B(v)\xi}{\lambda^{2/3}}\\
+ \big(|B(v)\eta|^2+F(v)\big)(\xi.\partial_{\eta})(\lambda^{-2/3})\Big].
\end{multline*}
Notice from (\ref{m2}) and (\ref{n5}) that
\begin{multline*}
\left|\frac{2B(v)\eta.B(v)\xi}{\lambda^{2/3}}\right| \leq \frac{2|B(v)\eta||B(v)\xi|}{\lambda^{2/3}} \leq \frac{2\sqrt{2}\lambda^{1/3}\lambda}{\lambda^{2/3}}  \leq 2\sqrt{2} \lambda^{2/3} \\ \leq 2\sqrt{2}
\big(|B(v)\eta|^2+F(v)\big),
\end{multline*}
on the support of the function
$$\psi'\left(\frac{|B(v)\eta|^2+F(v)}{\lambda^{2/3}}\right).$$
One can then deduce Lemma~\ref{lem6} from Lemma~\ref{lem5}.~$\Box$

\bigskip

\noindent
We notice from Lemma~\ref{lem5}, Lemma~\ref{lem6}, (\ref{n6}) and (\ref{n7}) that
$$\left|\{\xi.v,g\}-\frac{|B(v)\xi|^2}{\lambda^{4/3}}\psi\left(\frac{|B(v)\eta|^2+F(v)}{\lambda^{2/3}}\right)\right| \lesssim
1+|B(v)\eta|^2+F(v),$$
uniformly with respect to the parameter $\xi$ in $\rr^n$. It follows from (\ref{eq9}), (\ref{n8}) and the fact that the Wick quantization is a positive quantization (\ref{lay0.5}) that there exists $c_2>0$ such that for all $0<\eps \leq 1$, $\xi \in \rr^n$ and $u \in \mathcal{S}(\rr^n)$,
\begin{multline}\label{n9}
\frac{\eps}{4\pi}\left(\left[\frac{|B(v)\xi|^2}{\lambda^{4/3}}\psi\left(\frac{|B(v)\eta|^2+F(v)}{\lambda^{2/3}}\right)\right]^{\textrm{Wick}} u,u\right)
+\|B(v)\nabla_vu\|_{L^2}^2 +\big\|\sqrt{F(v)}u\big\|_{L^2}^2 \\ \leq C\|\langle B(v) \xi \rangle^{-s}Pu\|_{L^2}
\|\langle B(v) \xi \rangle^{s}u\|_{L^2}+c_2 \eps \big(\big[1+|B(v)\eta|^2+F(v)\big]^{\textrm{Wick}}u,u\big).
\end{multline}

\bigskip

\begin{lemma}\label{lem7}
There exists $c_{3}>0$ such that for all $\xi \in \rr^n$ and $u \in \mathcal{S}(\rr^n)$,
\begin{multline*}
\Big|\|B(v)\nabla_vu\|_{L^2}^2 +\big\|\sqrt{F(v)}u\big\|_{L^2}^2-\big(\big[4\pi^2|B(v)\eta|^2+F(v)\big]^{\emph{\textrm{Wick}}}u,u\big)\Big|  \\ \leq c_3\|\langle B(v) \xi \rangle^{-s}Pu\|_{L^2}
\|\langle B(v) \xi \rangle^{s}u\|_{L^2}.
\end{multline*}
\end{lemma}

\bigskip

\noindent
\textit{Proof of Lemma~\ref{lem7}.}
By using (\ref{lay1}) and (\ref{lay2}), we may write that
\begin{equation}\label{l2}
F(v)^{\textrm{Wick}}=F(v)+r_1^w,
\end{equation}
where $r_1$ is the real-valued symbol depending only on the variable $v$ given by
$$r_1(v)=\int_0^1\int_{\rr^{n}}{(1-\theta)\nabla_v^2\big(F(v)\big)(v+\theta \tilde{v}).\tilde{v}^2e^{-2\pi |\tilde{v}|^2}2^{n/2}d\tilde{v}d\theta}.$$
The notation
$$\nabla_v^2\big(F(v)\big)(v+\theta \tilde{v}),$$
denotes the second derivative with respect to the variable $v$ of the function $F(v)$ evaluated in the point $v+\theta \tilde{v} \in \rr^n$.
Since from (\ref{eq-3.1QQ}),
\begin{equation}\label{l2.5}
\big|\nabla_v^2\big(F(v)\big)\big| \lesssim \langle v \rangle^{\gamma} \textrm{ and }
\frac{\langle v\rangle}{\langle \theta \tilde{v} \rangle}\lesssim \langle v+\theta \tilde{v} \rangle \lesssim  \langle v\rangle \langle \theta \tilde{v} \rangle,
\end{equation}
we obtain by using again (\ref{eq-3.1QQ}) that
\begin{equation}\label{l3}
|(r_1^wu,u)|=|(r_1(v)u,u)|  \lesssim \|\langle v \rangle^{\gamma/2}u\|_{L^2}^2 \lesssim \big\|\sqrt{F(v)}u\big\|_{L^2}^2,
\end{equation}
because $|r_1(v)| \lesssim \langle v \rangle^{\gamma}$.
It follows from (\ref{eq1}), (\ref{l2}) and (\ref{l3}) that
\begin{equation}\label{l1}
\Big|\big\|\sqrt{F(v)}u\big\|_{L^2}^2-\big(F(v)^{\textrm{Wick}}u,u\big)\Big| \lesssim
\|\langle B(v) \xi \rangle^{-s}Pu\|_{L^2}
\|\langle B(v) \xi \rangle^{s}u\|_{L^2},
\end{equation}
uniformly with respect to the parameter $\xi \in \rr^n$.
By writing $D_v=i^{-1}\nabla_v$, we deduce from symbolic calculus and (\ref{eq-3QQ}) that
\begin{equation}\label{l4}
D_v.A(v)D_v=D_v.\big[A(v)\eta+ir_2(v)\big]^w,
\end{equation}
where $r_2$ is a real-valued symbol depending only on the variable $v$ and verifying
\begin{equation}\label{l5}
|r_2(v)| \lesssim \langle v \rangle^{\gamma+1}.
\end{equation}
It follows from the Cauchy-Schwarz inequality, (\ref{eq-3.1QQ}), (\ref{eq1}) and (\ref{l5}) that
\begin{align}
& \ \big|\big(D_v.r_2(v)u,u\big)\big|=\big|\big(\langle v \rangle^{-\gamma/2}r_2(v)u,\langle v \rangle^{\gamma/2}\nabla_v u\big)\big| \\
\lesssim & \  \|\langle v \rangle^{\frac{\gamma}{2}+1}u\|_{L^2}\|\langle v \rangle^{\gamma/2}\nabla_v u\|_{L^2}
\\ \lesssim & \  \|\langle v \rangle^{\frac{\gamma}{2}+1}u\|_{L^2}^2+\|\langle v \rangle^{\gamma/2}\nabla_v u\|_{L^2}^2
\lesssim \big\|\sqrt{F(v)}u\big\|_{L^2}^2 +\|\langle v \rangle^{\gamma/2}\nabla_v u\|_{L^2}^2
\\ \lesssim & \
\|\langle B(v) \xi \rangle^{-s}Pu\|_{L^2}
\|\langle B(v) \xi \rangle^{s}u\|_{L^2},\label{l6}
\end{align}
uniformly with respect to the parameter $\xi \in \rr^n$.
On the other hand, we may also write by using symbolic calculus, (\ref{eq-3QQ}) and (\ref{eq-2QQ}) that
\begin{equation}\label{l7}
D_v.\big(A(v)\eta\big)^w=\big[|B(v)\eta|^2+ir_3(v)\eta\big]^w=\big[|B(v)\eta|^2\big]^w+ir_3(v)D_v+r_4(v),
\end{equation}
where $r_3$ and $r_4$ are some real-valued symbols depending only on the variable $v$ and verifying
\begin{equation}\label{l8}
|r_3(v)| \lesssim \langle v \rangle^{\gamma+1} \textrm{ and }|r_4(v)| \lesssim \langle v \rangle^{\gamma}.
\end{equation}
It then follows from the Cauchy-Schwarz inequality, (\ref{eq-3.1QQ}), (\ref{eq1}) and (\ref{l8}) that
\begin{align}\label{l9}
& \ |(ir_3(v)D_vu,u)|=|(\langle v \rangle^{\gamma/2}\nabla_vu,\langle v \rangle^{-\gamma/2}r_3(v) u)| \\
\lesssim & \  \|\langle v \rangle^{\gamma/2}\nabla_v u\|_{L^2}\|\langle v \rangle^{\frac{\gamma}{2}+1}u\|_{L^2}
\lesssim   \|\langle v \rangle^{\gamma/2}\nabla_v u\|_{L^2}^2 +\|\langle v \rangle^{\frac{\gamma}{2}+1}u\|_{L^2}^2 \\
\lesssim & \ \|\langle v \rangle^{\gamma/2}\nabla_v u\|_{L^2}^2+ \big\|\sqrt{F(v)}u\big\|_{L^2}^2 \lesssim
\|\langle B(v) \xi \rangle^{-s}Pu\|_{L^2}
\|\langle B(v) \xi \rangle^{s}u\|_{L^2}
\end{align}
and
\begin{multline}\label{l10}
|(r_4(v)u,u)| \lesssim \|\langle v \rangle^{\gamma/2}u\|_{L^2}^2 \lesssim \big\|\sqrt{F(v)}u\big\|_{L^2}^2 \\ \lesssim \|\langle B(v) \xi \rangle^{-s}Pu\|_{L^2}
\|\langle B(v) \xi \rangle^{s}u\|_{L^2},
\end{multline}
uniformly with respect to the parameter $\xi \in \rr^n$.
One can therefore deduce from (\ref{eq-2QQ}), (\ref{l4}), (\ref{l6}), (\ref{l7}), (\ref{l9}) and (\ref{l10}) that
\begin{equation}\label{l11}
\big|\|B(v)\nabla_vu\|_{L^2}^2 -\big((|B(v)\eta|^2)^{w}u,u\big)\big|\lesssim  \|\langle B(v) \xi \rangle^{-s}Pu\|_{L^2}
\|\langle B(v) \xi \rangle^{s}u\|_{L^2},
\end{equation}
uniformly with respect to the parameter $\xi \in \rr^n$.
In view of (\ref{l1}) and (\ref{l11}), it remains to prove that
\begin{multline}\label{l12}
\big|4\pi^2\big((|B(v)\eta|^2)^{\textrm{Wick}}u,u\big)-\big((|B(v)\eta|^2)^{w}u,u\big)\big|\\ \lesssim \|\langle B(v) \xi \rangle^{-s}Pu\|_{L^2}
\|\langle B(v) \xi \rangle^{s}u\|_{L^2},
\end{multline}
uniformly with respect to the parameter $\xi \in \rr^n$;
in order to end the proof of Lemma~\ref{lem7}. Notice from (\ref{lay1}) and (\ref{lay2}) that we may write
\begin{equation}\label{l13}
(|B(v)\eta|^2)^{\textrm{Wick}}=\frac{1}{(2\pi)^2}(|B(v)\eta|^2)^{w}+R\Big(v,\frac{\eta}{2\pi}\Big)^w,
\end{equation}
with
$$R(v,\eta)=\int_0^1 \! \! \int_{\rr^{2n}}{(1-\theta)\nabla_{v,\eta}^2\big(|B(v)\eta|^2\big)(v+\theta \tilde{v},\eta+\theta \tilde{\eta}).(\tilde{v},\tilde{\eta})^2e^{-2\pi |(\tilde{v},\tilde{\eta})|^2}2^{n}
d\tilde{v}d\tilde{\eta}d\theta}.$$
The factor $2\pi$ comes from the fact that we are using here the normalization of the Weyl quantization defined in (\ref{weyl1}) which differs from the one used in the appendix (See (\ref{lay3})).
Define
$$R_1=\int_0^1 \! \! \! \int_{\rr^{2n}}{(1-\theta)\nabla_{\eta}^2\big(|B(v)\eta|^2\big)(v+\theta \tilde{v},\eta+\theta \tilde{\eta}).\tilde{\eta}^2e^{-2\pi |(\tilde{v},\tilde{\eta})|^2}2^{n}
d\tilde{v}d\tilde{\eta}d\theta},$$
$$R_2=\int_0^1 \! \! \!\int_{\rr^{2n}}{(1-\theta)\nabla_{v}\nabla_{\eta}\big(|B(v)\eta|^2\big)(v+\theta \tilde{v},\eta+\theta \tilde{\eta}).(\tilde{v},\tilde{\eta})^2e^{-2\pi |(\tilde{v},\tilde{\eta})|^2}2^{n}
d\tilde{v}d\tilde{\eta}d\theta}$$
and
$$R_3=\int_0^1 \! \! \!\int_{\rr^{2n}}{(1-\theta)\nabla_{v}^2\big(|B(v)\eta|^2\big)(v+\theta \tilde{v},\eta+\theta \tilde{\eta}).\tilde{v}^2e^{-2\pi |(\tilde{v},\tilde{\eta})|^2}2^{n}
d\tilde{v}d\tilde{\eta}d\theta}.$$
We first notice that the symbol $R_1$ only depends on the variable $v$. We then notice from (\ref{eq-1QQ}) that
\begin{equation}\label{l14}
|\nabla_{\eta}^2\big(|B(v)\eta|^2\big)| \lesssim \langle v \rangle^{\gamma+2}.
\end{equation}
By using (\ref{l2.5}) as above, one can estimate the function $R_1(v)$  from above as
$$|R_1(v)| \lesssim \langle v \rangle^{\gamma+2}.$$
This therefore implies that
\begin{multline}\label{l14.1}
\Big|\Big(R_1\Big(v,\frac{\eta}{2\pi}\Big)^wu,u\Big)\Big| =|(R_1(v)u,u)| \lesssim \|\langle v \rangle^{\frac{\gamma}{2}+1}u\|_{L^2}^2  \lesssim \big\|\sqrt{F(v)}u\big\|_{L^2}^2
 \\ \lesssim \|\langle B(v) \xi \rangle^{-s}Pu\|_{L^2}
\|\langle B(v) \xi \rangle^{s}u\|_{L^2},
\end{multline}
uniformly with respect to the parameter $\xi \in \rr^n$;
by using (\ref{eq-3.1QQ}) and (\ref{eq1}). We then notice always from (\ref{eq-1QQ}) that
$$|\nabla_v\nabla_{\eta}\big(|B(v)\eta|^2\big)| \lesssim \langle v \rangle^{\gamma+1}|\eta|,$$
and that the term $\nabla_v\nabla_{\eta}\big(|B(v)\eta|^2\big)$ is linear in the variable $\eta$. By using again (\ref{eq-1QQ}), this shows that we may write the symbol $R_2\big(v,\eta/(2\pi)\big)$ as
$$R_2\Big(v,\frac{\eta}{2\pi}\Big)=a(v)\eta+b(v),$$
where $a$ and $b$ are smooth functions verifying for all $\alpha, \beta \in \nn^n$,
$$|\partial_v^{\alpha}a(v)|+|\partial_v^{\beta}b(v)| \lesssim \langle v \rangle^{\gamma+1}.$$
Thus, we may write by using symbolic calculus that
$$R_2\Big(v,\frac{\eta}{2\pi}\Big)^w=a(v)D_v+\tilde{b}(v),$$
where $\tilde{b}$ is a smooth function fulfilling for all $\alpha \in \nn^n$,
$$|\partial_v^{\alpha}\tilde{b}(v)| \lesssim \langle v \rangle^{\gamma+1}.$$
It follows from the Cauchy-Schwarz inequality, (\ref{eq-3.1QQ}) and (\ref{eq1}) that
\begin{align}\label{l15}
& \ \Big|\Big(R_2\Big(v,\frac{\eta}{2\pi}\Big)^w u,u\Big)\Big| \lesssim \|\langle v \rangle^{\gamma/2}\nabla_v u\|_{L^2} \|\langle v \rangle^{-\gamma/2}a(v)u\|_{L^2}
+\|\langle v \rangle^{\frac{\gamma+1}{2}} u\|_{L^2}^2 \\
\lesssim  &\
\|\langle v \rangle^{\gamma/2}\nabla_v u\|_{L^2}^2+\|\langle v \rangle^{\frac{\gamma}{2}+1} u\|_{L^2}^2
\lesssim  \|\langle v \rangle^{\gamma/2}\nabla_v u\|_{L^2}^2+\big\|\sqrt{F(v)} u\big\|_{L^2}^2
\\ \lesssim & \
\|\langle B(v) \xi \rangle^{-s}Pu\|_{L^2}
\|\langle B(v) \xi \rangle^{s}u\|_{L^2},
\end{align}
uniformly with respect to the parameter $\xi \in \rr^n$.
We finally notice from (\ref{eq-1QQ}) that
$$ |\nabla_{v}^2\big(|B(v)\eta|^2\big)| \lesssim \langle v \rangle^{\gamma}|\eta|^2,$$
and that the term $\nabla_{v}^2\big(|B(v)\eta|^2\big)$ is quadratic in the variable $\eta$. By using again (\ref{eq-1QQ}), this shows that we may write the symbol
$R_3\big(v,\eta/(2\pi)\big)$ as
$$R_3\Big(v,\frac{\eta}{2\pi}\Big)=a(v)|\eta|^2+b(v)\eta+c(v),$$
where $a$, $b$ and $c$ are smooth functions verifying for all $\alpha_1, \alpha_2, \alpha_3 \in \nn^n$,
$$|\partial_v^{\alpha_1}a(v)|+|\partial_v^{\alpha_2}b(v)|+|\partial_v^{\alpha_3}c(v)| \lesssim \langle v \rangle^{\gamma}.$$
Thus, we may write by using symbolic calculus that
$$R_3\Big(v,\frac{\eta}{2\pi}\Big)^w=\nabla_v.\tilde{a}(v)\nabla_v+\tilde{b}(v)\nabla_v+\tilde{c}(v),$$
where $\tilde{a}$, $\tilde{b}$, $\tilde{c}$ are smooth functions verifying for all $\alpha_1, \alpha_2, \alpha_3 \in \nn^n$,
$$|\partial_v^{\alpha_1}\tilde{a}(v)|+|\partial_v^{\alpha_2}\tilde{b}(v)|+|\partial_v^{\alpha_3}\tilde{c}(v)| \lesssim \langle v \rangle^{\gamma}.$$
It follows from the Cauchy-Schwarz inequality, (\ref{eq-3.1QQ}) and (\ref{eq1}) that
\begin{align}\label{l16}
& \ \Big|\Big(R_3\Big(v,\frac{\eta}{2\pi}\Big)^w u,u\Big)\Big|\\
\lesssim & \ \|\langle v \rangle^{\gamma/2}\nabla_v u\|_{L^2}\|\langle v \rangle^{-\gamma}\tilde{a}(v)\langle v \rangle^{\gamma/2}\nabla_v u\|_{L^2}
+\|\langle v \rangle^{\gamma/2} u\|_{L^2}^2 \\
+& \ \|\langle v \rangle^{\gamma/2}\nabla_v u\|_{L^2}\|\langle v \rangle^{-\gamma/2}\tilde{b}(v) u\|_{L^2}
\lesssim \|\langle v \rangle^{\gamma/2}\nabla_v u\|_{L^2}^2+\|\langle v \rangle^{\frac{\gamma}{2}+1} u\|_{L^2}^2 \\
\lesssim & \
\|\langle v \rangle^{\gamma/2}\nabla_v u\|_{L^2}^2+\big\|\sqrt{F(v)} u\big\|_{L^2}^2\lesssim
\|\langle B(v) \xi \rangle^{-s}Pu\|_{L^2}
\|\langle B(v) \xi \rangle^{s}u\|_{L^2},
\end{align}
uniformly with respect to the parameter $\xi \in \rr^n$.
One can then deduce from (\ref{l14.1}), (\ref{l15}) and (\ref{l16}) that
$$\Big|\Big(R\Big(v,\frac{\eta}{2\pi}\Big)^w u,u\Big)\Big| \lesssim \|\langle B(v) \xi \rangle^{-s}Pu\|_{L^2}
\|\langle B(v) \xi \rangle^{s}u\|_{L^2},$$
uniformly with respect to the parameter $\xi \in \rr^n$.
According to (\ref{l13}), this proves the estimate (\ref{l12}) and ends the proof of Lemma~\ref{lem7}.~$\Box$

\subsubsection{Proof of Proposition \ref{prop1}}

We now take advantage of all the results
proved previously. In particular,  it follows from Lemma~\ref{lem7} and (\ref{n9}) that there exist $0<\eps_0 \leq 1$ and $c_4>0$ such that for all $\xi \in \rr^n$, $0<\eps \leq \eps_0$ and $u \in \mathcal{S}(\rr^n)$,
\begin{multline*}
\eps\left(\left[\frac{|B(v)\xi|^2}{\lambda^{4/3}}\psi\left(\frac{|B(v)\eta|^2+F(v)}{\lambda^{2/3}}\right)\right]^{\textrm{Wick}} u,u\right)
+\|B(v)\nabla_vu\|_{L^2}^2\\ +\big\|\sqrt{F(v)}u\big\|_{L^2}^2  \leq c_4 \|\langle B(v) \xi \rangle^{-s}Pu\|_{L^2}
\|\langle B(v) \xi \rangle^{s}u\|_{L^2}+c_4\|u\|_{L^2}^2.
\end{multline*}
Another use of Lemma~\ref{lem7} shows that there exists $c_5>0$ such that for all $\xi \in \rr^n$, $0<\eps \leq \eps_0$ and $u \in \mathcal{S}(\rr^n)$,
\begin{multline}\label{n10}
\eps\left(\left[\frac{|B(v)\xi|^2}{\lambda^{4/3}}\psi\left(\frac{|B(v)\eta|^2+F(v)}{\lambda^{2/3}}\right)\right]^{\textrm{Wick}} u,u\right)+\|u\|_{L^2}^2
\\ +\big(\big[4\pi^2|B(v)\eta|^2+F(v)\big]^{\textrm{Wick}}u,u\big)  \leq c_5 \|\langle B(v) \xi \rangle^{-s}Pu\|_{L^2}
\|\langle B(v) \xi \rangle^{s}u\|_{L^2}+c_5\|u\|_{L^2}^2.
\end{multline}
Notice from (\ref{m2}) and (\ref{n5}) that
\begin{align*}
& \ \eps\frac{|B(v)\xi|^2}{\lambda^{4/3}}\psi\left(\frac{|B(v)\eta|^2+F(v)}{\lambda^{2/3}}\right)+4\pi^2|B(v)\eta|^2+F(v)+1 \gtrsim \\
& \ \eps\frac{|B(v)\xi|^2+|B(v)\eta|^2+F(v)+1}{\lambda^{4/3}}\psi\left(\frac{|B(v)\eta|^2+F(v)}{\lambda^{2/3}}\right)\\
 & \ +\big(|B(v)\eta|^2+F(v)+1\big)\left[1-\psi\left(\frac{|B(v)\eta|^2+F(v)}{\lambda^{2/3}}\right)\right] \geq \\
& \ \eps\frac{\lambda^2}{\lambda^{4/3}}\psi\left(\frac{|B(v)\eta|^2+F(v)}{\lambda^{2/3}}\right)
 +\lambda^{2/3}\left[1-\psi\left(\frac{|B(v)\eta|^2+F(v)}{\lambda^{2/3}}\right)\right] \geq \eps \lambda^{2/3},
\end{align*}
when $0<\eps \leq1$; since
$$|B(v)\eta|^2+F(v) \geq \lambda^{2/3},$$
on the support of the function
$$1-\psi\left(\frac{|B(v)\eta|^2+F(v)}{\lambda^{2/3}}\right).$$
This is the crucial step where we use that the multiplier (\ref{m1}) creates the good term
$$\frac{|B(v)\xi|^2}{\lambda^{4/3}}\psi\left(\frac{|B(v)\eta|^2+F(v)}{\lambda^{2/3}}\right),$$
in order to control the quantity $\lambda^{2/3}$ in the region of the phase space where
$$|B(v)\eta|^2+F(v) \leq \lambda^{2/3}.$$
By using again that the Wick quantization is a positive quantization (\ref{lay0.5}), we deduce that there exists $c_6>0$ such that for all $\xi \in \rr^n$ and $u \in \mathcal{S}(\rr^n)$,
\begin{equation}\label{n11}
\big((\lambda^{2/3})^{\textrm{Wick}}u,u\big)  \leq c_6 \|\langle B(v) \xi \rangle^{-s}Pu\|_{L^2}
\|\langle B(v) \xi \rangle^{s}u\|_{L^2}+c_6\|u\|_{L^2}^2.
\end{equation}
Since from (\ref{m2}),
$$\langle B(v)\xi \rangle^{2/3} \lesssim \lambda^{2/3},$$
it follows from (\ref{lay0.5}) that there exists $c_7>0$ such that for all $\xi \in \rr^n$ and $u \in \mathcal{S}(\rr^n)$,
\begin{equation}\label{n12}
\big(\big(\langle B(v)\xi\rangle^{2/3}\big)^{\textrm{Wick}}u,u\big) \leq c_7 \|\langle B(v) \xi \rangle^{-s}Pu\|_{L^2}
\|\langle B(v) \xi \rangle^{s}u\|_{L^2}+c_7\|u\|_{L^2}^2.
\end{equation}
Notice from (\ref{lay1bis}) and (\ref{lay2bis}) that we have
$$m(v,\xi)^2=\big(\langle B(v)\xi\rangle^{2/3}\big)^{\textrm{Wick}},$$
where we recall that the quantity $m(v,\xi)$ is defined in the statement of Proposition~\ref{prop1}.
We finally obtain from (\ref{eq1}) and (\ref{n12}) that  there exists $c_8>0$ such that for all $\xi \in \rr^n$ and $u \in \mathcal{S}(\rr^n)$,
\begin{multline*}
\|m(v,\xi)u\|_{L^2}^2 +\|B(v)\nabla_vu\|_{L^2}^2+\big\|\sqrt{F(v)}u\big\|_{L^2}^2\\ \leq c_8 \|\langle B(v) \xi \rangle^{-s}Pu\|_{L^2}
\|\langle B(v) \xi \rangle^{s}u\|_{L^2}+c_8\|u\|_{L^2}^2.
\end{multline*}
This ends the proof of Proposition~\ref{prop1}. $\Box$

\bigskip

\noindent
The next proposition follows directly from Proposition~\ref{prop1} and gives first non-optimal hypoelliptic estimates fulfilled by generalized linear Landau-type operators.  Notice that this estimate will be instrumental in the proof of Theorem~\ref{mainlandau}.

\bigskip

\begin{proposition}\label{prop2}
Let $P$ be a generalized linear Landau-type operator fulfilling the assumptions \emph{(\ref{eq-3QQ})}, \emph{(\ref{eq-3.1QQ})}, \emph{(\ref{eq-2QQ})}, \emph{(\ref{eq-1QQ})} and \emph{(\ref{eq0QQ})}. Then, there exists a constant $C>0$ such that for all $\xi \in \rr^n$ and $u \in \mathcal{S}(\rr_{v}^{n})$,
$$\|m(v,\xi)u\|_{L^2}^2 +\|B(v)\nabla_vu\|_{L^2}^2+\big\|\sqrt{F(v)}u\big\|_{L^2}^2\leq C\big(\|Pu\|_{L^2}^2+\|u\|_{L^2}^2\big)$$
and
$$\|\langle v \rangle^{\gamma/6}|\xi|^{1/3}u\|_{L^2}^2 +\|B(v)\nabla_vu\|_{L^2}^2+\|\langle v \rangle^{\frac{\gamma}{2}+1}u\|_{L^2}^2\leq C\big(\|Pu\|_{L^2}^2+\|u\|_{L^2}^2\big),$$
with
$$m(v,\xi)=\Big(\int_{\rr^n}\langle B(v+\tilde{v})\xi \rangle^{2/3}e^{-2\pi\tilde{v}^2}2^{n/2}d\tilde{v}\Big)^{1/2},$$
where the notation $\|\cdot\|$ stands for the $L^2(\rr_{v}^{n})$-norm.
\end{proposition}

\bigskip

\noindent
\textit{Proof of Proposition~\ref{prop2}.} We deduce directly the first estimate from Proposition~\ref{prop1} by taking $s=0$. To prove the second estimate, we just need to use (\ref{eq-3.1QQ}) and to prove that
$$\langle v \rangle^{\gamma/3}|\xi|^{2/3} \lesssim m(v,\xi)^2.$$
By using that
\begin{equation}\label{li3}
\frac{\langle a \rangle}{\langle b \rangle} \lesssim \langle a+b \rangle \lesssim \langle a \rangle \langle b \rangle,
\end{equation}
it follows from (\ref{eq0QQ}) that
\begin{multline*}
m(v,\xi)^2 \gtrsim \int_{\rr^n}|B(v+\tilde{v})\xi|^{2/3}e^{-2\pi\tilde{v}^2}d\tilde{v} \gtrsim \int_{\rr^n}\langle v+\tilde{v} \rangle^{\gamma/3}
|\xi|^{2/3}e^{-2\pi\tilde{v}^2}d\tilde{v} \\ \gtrsim \int_{\rr^n}\frac{\langle v \rangle^{\gamma/3}}{\langle \tilde{v} \rangle^{|\gamma|/3}}
|\xi|^{2/3}e^{-2\pi\tilde{v}^2}d\tilde{v} \gtrsim \langle v \rangle^{\gamma/3}|\xi|^{2/3}.
\end{multline*}
This ends the proof of Proposition~\ref{prop2}.~$\Box$

\subsection{Hypoelliptic estimates for linear Landau-type operators}\label{la2}
This section is devoted to the proof of Theorem \ref{mainlandau}.
We consider a linear Landau-type operator
\begin{equation}\label{fili11}
P=iv.D_x+ D_v.\lambda(v)D_{v}+Q.\mu(v) Q+F(v), \ x,v \in \rr^{3};
\end{equation}
with $D_x=i^{-1}\partial_x$, $D_v=i^{-1}\partial_v$, $\gamma \in [-3,1]$, where $Q = v \wedge D_v$ stands also for the vector-valued operator defined by the Weyl quantization of the vector-valued symbol $v \wedge \eta$; and where $F$, $\lambda$ and $\mu$ are some positive functions satisfying (\ref{eq-3.1Q}) and (\ref{w2Q}).
As mentioned previously, a linear Landau-type operator is a generalized linear Landau-type operator with $B(v)$ explicitly defined in (\ref{an2QQ}).

Starting from the a priori estimate proved in the previous section, we shall now establish sharp hypoelliptic estimates with loss of $4/3$ derivatives for linear Landau-type operators. As in the Fokker-Planck case, we split the proof of Theorem \ref{mainlandau} into two parts deriving separately estimates for the spatial derivatives and the velocity variables.

\subsubsection{\underline{Spatial derivatives estimates}}

The aim of this subsection is to give a proof of the following Proposition:

\bigskip

\begin{proposition}\label{propla}
Let $P$ be the linear Landau-type operator defined in \emph{(\ref{fili11})}. Then, there exists $C>0$ such that for all $u \in \mathcal{S}(\rr_{x,v}^{6})$,
$$\|\langle B(v)D_x \rangle^{2/3}u\|_{L^2}^2 \leq C( \|Pu\|_{L^2}^2+\|u\|_{L^2}^2)$$
and
$$\|\langle v \rangle^{\gamma/3}|D_x|^{2/3}u\|_{L^2}^2 +\|\langle v \rangle^{\gamma/3}|v \wedge D_x|^{2/3}u\|_{L^2}^2  \leq C(
 \|Pu\|_{L^2}^2+\|u\|_{L^2}^2),$$
where $\|\cdot\|_{L^2}$ stands for the $L^2(\rr_{x,v}^{6})$-norm.
\end{proposition}

\bigskip

\noindent
We obtain from Proposition~\ref{prop2} and (\ref{eq-3.1QQ}) the estimate
\begin{equation}\label{li4b}
\|m(v,\xi)u\|_{L^2}^2 +\|B(v)\nabla_vu\|_{L^2}^2+\|\langle v \rangle^{\frac{\gamma}{2}+1}u\|_{L^2}^2 \lesssim  \|Pu\|_{L^2}^2+\|u\|_{L^2}^2,
\end{equation}
uniformly with respect to the parameter $\xi \in \rr^n$; where the notation $\|\cdot\|_{L^2}$ stands for the $L^2(\rr_{v}^{n})$-norm and
$$m(v,\xi)=\Big(\int_{\rr^3}\langle B(v+\tilde{v})\xi \rangle^{2/3}e^{-2\pi\tilde{v}^2}2^{3/2}d\tilde{v}\Big)^{1/2}.$$
In the specific case of a linear Landau-type operator, one can simply estimate from below the term $m(v,\xi)$.

\bigskip

\begin{lemma}\label{lem8}There exists $C>0$ such that for all $v \in \rr^3$ and $\xi \in \rr^3$,
$$\langle B(v)\xi\rangle^{1/3} \leq C m(v,\xi).$$
\end{lemma}

\bigskip

\noindent
\textit{Proof of Lemma~\ref{lem8}.} According to (\ref{rih1}), we may write that
$$m(v,\xi)^2=\int_{\rr^3}\big(1+|\sqrt{\lambda(v+\tilde{v})}\xi|^2+
|\sqrt{\mu(v+\tilde{v})} (v+\tilde{v}) \wedge \xi|^2 \big)^{1/3}e^{-2\pi\tilde{v}^2}2^{3/2}d\tilde{v}.$$
This implies that
\begin{multline}\label{li1}
m(v,\xi)^2 \gtrsim 1+\int_{\rr^3}|\sqrt{\mu(v+\tilde{v})} (v+\tilde{v}) \wedge \xi|^{2/3}e^{-2\pi\tilde{v}^2}d\tilde{v}
\\ +\int_{\rr^3}|\sqrt{\lambda(v+\tilde{v})}\xi|^{2/3}e^{-2\pi\tilde{v}^2}d\tilde{v}. \end{multline}
By using (\ref{li3}),
one can then notice from (\ref{w2Q}) that
\begin{multline*}
\int_{\rr^3}|\sqrt{\lambda(v+\tilde{v})}\xi|^{2/3}e^{-2\pi\tilde{v}^2}d\tilde{v} \gtrsim
\int_{\rr^3}\langle v+\tilde{v}\rangle^{\gamma/3}|\xi|^{2/3}e^{-2\pi\tilde{v}^2}d\tilde{v} \\  \gtrsim \int_{\rr^3}\frac{\langle v \rangle^{\gamma/3}}{\langle \tilde{v} \rangle^{|\gamma|/3}}|\xi|^{2/3}e^{-2\pi\tilde{v}^2}d\tilde{v} \gtrsim \langle v \rangle^{\gamma/3}|\xi|^{2/3}.
\end{multline*}
We therefore obtain that
\begin{equation}\label{li2}
m(v,\xi)^2 \gtrsim 1+\langle v \rangle^{\gamma/3}|\xi|^{2/3}+\int_{\rr^3}|\sqrt{\mu(v+\tilde{v})} (v+\tilde{v}) \wedge \xi|^{2/3}e^{-2\pi\tilde{v}^2}d\tilde{v}.
\end{equation}
On the other hand, we also deduce from (\ref{w2Q}) and (\ref{li3}) that
\begin{multline*}
\int_{\rr^3}|\sqrt{\mu(v+\tilde{v})} (v+\tilde{v}) \wedge \xi|^{2/3}e^{-2\pi\tilde{v}^2}d\tilde{v} \gtrsim
\int_{\rr^3}\langle v+\tilde{v} \rangle^{\gamma/3} |(v+\tilde{v}) \wedge \xi|^{2/3}e^{-2\pi\tilde{v}^2}d\tilde{v}\\
\gtrsim \int_{\overline{B(0,1)}}\langle v+\tilde{v} \rangle^{\gamma/3} |(v+\tilde{v}) \wedge \xi|^{2/3}d\tilde{v}
\gtrsim  \int_{\overline{B(0,1)}}\frac{\langle v\rangle^{\gamma/3}}{\langle \tilde{v}\rangle^{|\gamma|/3}} |(v+\tilde{v}) \wedge \xi|^{2/3}d\tilde{v},
\end{multline*}
where $\overline{B(0,1)}$ stands for the closed unit ball in $\rr^3$. By noticing that we have
$$ |(v+\tilde{v}) \wedge \xi| \geq  |v \wedge \xi|- |\tilde{v} \wedge \xi| \geq  |v \wedge \xi|- |\tilde{v}| |\xi| \geq |v \wedge \xi|- |\xi|\geq \frac{1}{2}|v \wedge \xi|,$$
when $|\tilde{v}| \leq 1$ and $2|\xi| \leq |v \wedge \xi|$, it follows that
$$\int_{\rr^3}|\sqrt{\mu(v+\tilde{v})} (v+\tilde{v}) \wedge \xi|^{2/3}e^{-2\pi\tilde{v}^2}d\tilde{v} \gtrsim
\langle v\rangle^{\gamma/3}|v\wedge \xi|^{2/3},$$
when $2|\xi| \leq |v \wedge \xi|$. Since
$$\langle v \rangle^{\gamma/3}|\xi|^{2/3} \gtrsim \langle v\rangle^{\gamma/3}|v\wedge \xi|^{2/3},$$
when $2|\xi| \geq |v \wedge \xi|$, it follows from (\ref{eq-3.1Q}), (\ref{w2Q}), (\ref{rih1}) and (\ref{li2}) that
\begin{multline}\label{li4}
m(v,\xi)^2 \gtrsim 1+\langle v \rangle^{\gamma/3}|\xi|^{2/3}+\langle v\rangle^{\gamma/3}|v\wedge \xi|^{2/3}
\gtrsim 1+|\sqrt{\lambda(v)}\xi|^{2/3}\\ +|\sqrt{\mu(v)}v \wedge \xi|^{2/3}
\gtrsim \big(1+|\sqrt{\lambda(v)}\xi|^{2}+|\sqrt{\mu(v)}v \wedge \xi|^{2}\big)^{1/3}
=\langle B(v)\xi \rangle^{2/3}.
\end{multline}
This proves Lemma~\ref{lem8}.~$\Box$

\bigskip

\noindent
By coming back to the direct side in the $x$ variable, we deduce from Lemma~\ref{lem8} and (\ref{li4b}) that there exists $C>0$ such that for all $u \in \mathcal{S}(\rr_{x,v}^{6})$,
\begin{equation}\label{and1}
\|\langle B(v)D_x \rangle^{1/3}u\|_{L^2}^2 +\|B(v)\nabla_vu\|_{L^2}^2+\|\langle v \rangle^{\frac{\gamma}{2}+1}u\|_{L^2}^2\leq C(\|Pu\|_{L^2}^2+\|u\|_{L^2}^2),
\end{equation}
where $\|\cdot\|$ stands for the $L^2(\rr_{x,v}^{6})$-norm. According to (\ref{w2Q}), (\ref{eq0QQ}) and (\ref{li4}),
this implies in particular that
\begin{multline}\label{and2}
\|\langle v \rangle^{\gamma/6}|D_x|^{1/3}u\|_{L^2}^2 +\|\langle v \rangle^{\gamma/6}|v \wedge D_x|^{1/3}u\|_{L^2}^2+\|\langle v\rangle^{\gamma/2}\nabla_vu\|_{L^2}^2\\ +\|\langle v \rangle^{\frac{\gamma}{2}+1}u\|_{L^2}^2  \lesssim  \|Pu\|_{L^2}^2+\|u\|_{L^2}^2.
\end{multline}
We shall now improve these estimates and prove Proposition~\ref{propla} by using an argument of commutation.

\bigskip

\begin{lemma}\label{lemme1}
For any $s \in \rr$, we have
$$\|\langle B(v)\xi \rangle^{-s}[D_v.\lambda(v)D_v,\langle B(v)\xi \rangle^{s}]u\|_{L^2}^2 \lesssim \|\langle v \rangle^{\gamma}\nabla_vu\|_{L^2}^2+\|\langle v \rangle^{\gamma}u\|_{L^2}^2,$$
uniformly with respect to the parameter $\xi$ in $\rr^3$; where $\|\cdot\|_{L^2}$ stands for the $L^2(\rr_{v}^{3})$ norm.
\end{lemma}

\bigskip

\noindent
\textit{Proof of Lemma~\ref{lemme1}.}
We may write that
\begin{multline}\label{dio1}
\langle B(v)\xi \rangle^{-s}[D_v.\lambda(v)D_v,\langle B(v)\xi \rangle^{s}]=\langle B(v)\xi \rangle^{-s}\Big(D_v.\lambda(v)[D_v,\langle B(v)\xi \rangle^{s}]\\ +
[D_v,\langle B(v)\xi \rangle^{s}].\lambda(v)D_v\Big).
\end{multline}
Symbolic calculus shows that
\begin{equation}\label{dio3}
[D_v,\langle B(v)\xi \rangle^{s}]=\frac{1}{i}\nabla_{v}\big(\langle B(v)\xi \rangle^{s}\big).
\end{equation}
It follows from Lemma~\ref{lem0} and (\ref{eq-3.1Q}) that
\begin{multline}\label{dio2}
\|\langle B(v)\xi \rangle^{-s}[D_v,\langle B(v)\xi \rangle^{s}].\lambda(v)D_vu\|_{L^2}^2\\ =\|\langle B(v)\xi \rangle^{-s}\nabla_{v}\big(\langle B(v)\xi \rangle^{s}\big).\lambda(v)D_vu\|_{L^2}^2   \lesssim \|\langle v \rangle^{\gamma}\nabla_v u\|_{L^2}^2,
\end{multline}
since
$$\big|\langle B(v)\xi \rangle^{-s}\nabla_{v}\big(\langle B(v)\xi \rangle^{s}\big)\big| \lesssim 1,$$
uniformly with respect to the parameter $\xi$ in $\rr^3$.
Notice from Lemma~\ref{lem0} and (\ref{eq-3.1Q}) that
\begin{equation}\label{dio3b}
\lambda(v)\nabla_{v}\big(\langle B(v)\xi \rangle^{s}\big) \in S\big(\langle v \rangle^{\gamma}\langle B(v)\xi \rangle^{s},dv^2+d\eta^2\big).
\end{equation}
Keeping in mind (\ref{dio3}), another use of symbolic calculus shows that we may write
\begin{multline}\label{dio3.5}
\langle B(v)\xi \rangle^{-s}D_v.\lambda(v)[D_v,\langle B(v)\xi \rangle^{s}]=\frac{1}{i}
\langle B(v)\xi \rangle^{-s}\lambda(v)\nabla_{v}\big(\langle B(v)\xi \rangle^{s}\big).D_v\\ +b(v,\xi),
\end{multline}
where $b$ is a smooth function depending only on the variable $v$ and the parameter $\xi \in \rr^3$, and verifying
\begin{equation}\label{dio4}
|b(v,\xi)| \lesssim  \langle v \rangle^{\gamma},
\end{equation}
uniformly with respect to the parameter $\xi$ in $\rr^3$.
It follows from (\ref{dio3b}), (\ref{dio3.5}) and (\ref{dio4}) that
\begin{multline}\label{dio5}
\|\langle B(v)\xi \rangle^{-s}D_v.\lambda(v)[D_v,\langle B(v)\xi \rangle^{s}]u\|_{L^2}^2 \\ \lesssim \|\langle B(v)\xi \rangle^{-s}\lambda(v)\nabla_{v}\big(\langle B(v)\xi \rangle^{s}\big).D_vu\|_{L^2}^2 +\|b(v,\xi)u\|_{L^2}^2
\lesssim \|\langle v \rangle^{\gamma}\nabla_vu\|_{L^2}^2+\|\langle v \rangle^{\gamma}u\|_{L^2}^2,
\end{multline}
since
$$\big|\langle B(v)\xi \rangle^{-s}\lambda(v)\nabla_{v}\big(\langle B(v)\xi \rangle^{s}\big)\big| \lesssim \langle v \rangle^{\gamma},$$
uniformly with respect to the parameter $\xi$ in $\rr^3$.
One can then deduce from (\ref{dio1}), (\ref{dio2}) and (\ref{dio5}) the estimate of Lemma~\ref{lemme1}.~$\Box$

\bigskip

\begin{lemma}\label{lemme2}
For any $s \in \rr$, we have
$$\big\|\langle B(v)\xi \rangle^{-s}\big[Q.\mu(v)Q,\langle B(v)\xi \rangle^{s}\big]u\big\|_{L^2}^2 \lesssim
\|\langle v \rangle^{\gamma+1}Qu\|_{L^2}^2+\|\langle v \rangle^{\gamma+2}u\|_{L^2}^2,$$
uniformly with respect to the parameter $\xi$ in $\rr^3$; where $\|\cdot\|_{L^2}$ stands for the $L^2(\rr_{v}^{3})$ norm.

\end{lemma}

\bigskip

\noindent
\textit{Proof of Lemma~\ref{lemme2}.}
We may write that
\begin{multline}\label{dio1b}
\langle B(v)\xi \rangle^{-s}\big[Q.\mu(v) Q,\langle B(v)\xi \rangle^{s}\big]=\langle B(v)\xi \rangle^{-s}
\big(Q.\mu(v)[Q,\langle B(v)\xi \rangle^{s}]\\ +
[Q,\langle B(v)\xi \rangle^{s}].\mu(v)Q\big).
\end{multline}
Since
\begin{equation}\label{ed0}
Q=(v \wedge \eta)^w=\left(
\begin{array}{c}
v_2D_{v_3}-v_3D_{v_2}\\
v_3D_{v_1}-v_1D_{v_3}\\
v_1D_{v_2}-v_2D_{v_1}
\end{array}
\right),
\end{equation}
symbolic calculus and Lemma~\ref{lem0} show that
\begin{equation}\label{ed1}
[Q,\langle B(v)\xi \rangle^s]=a(v,\xi),
\end{equation}
where $a$ is a smooth function depending only on the variable $v$ and the parameter $\xi$ in $\rr^3$, and verifying
\begin{equation}\label{ed2}
a \in S\big(\langle v \rangle \langle B(v)\xi \rangle^s,dv^2+d\eta^2\big),
\end{equation}
uniformly with respect to the parameter $\xi$ in $\rr^3$.
It follows from (\ref{eq-3.1Q}), (\ref{ed1}) and (\ref{ed2}) that
\begin{equation}\label{ed3}
\|\langle B(v)\xi \rangle^{-s}[Q,\langle B(v)\xi \rangle^{s}].\mu(v)Qu\|_{L^2}^2 \lesssim
\|\langle v \rangle^{\gamma+1}Qu\|_{L^2}^2,
\end{equation}
since
$$|\langle B(v)\xi \rangle^{-s}\mu(v)a(v,\xi)| \lesssim \langle v \rangle^{\gamma+1},$$
uniformly with respect to the parameter $\xi$ in $\rr^3$.
Setting
\begin{equation}\label{ed-11}
b(v,\xi)=\mu(v)[Q,\langle B(v)\xi \rangle^{s}],
\end{equation}
it follows from (\ref{eq-3.1Q}), (\ref{ed1}) and (\ref{ed2}) that
\begin{equation}\label{ed4}
b \in S\big(\langle v \rangle^{\gamma+1} \langle B(v)\xi \rangle^s,dv^2+d\eta^2\big),
\end{equation}
uniformly with respect to the parameter $\xi$ in $\rr^3$. According to (\ref{ed0}), (\ref{ed-11}) and (\ref{ed4}), symbolic calculus shows that
$$\big[Q,\mu(v)[Q,\langle B(v)\xi \rangle^{s}]\big]=c(v,\xi),$$
where $c$ is a smooth function depending only on the variable $v$ and the parameter $\xi$ in $\rr^3$, and verifying
\begin{equation}\label{ed5}
c \in S\big(\langle v \rangle^{\gamma+2} \langle B(v)\xi \rangle^s,dv^2+d\eta^2\big),
\end{equation}
uniformly with respect to the parameter $\xi$ in $\rr^3$. This implies that
\begin{multline}\label{ed7}
\big\|\langle B(v)\xi \rangle^{-s}\big[Q,\mu(v)[Q,\langle B(v)\xi \rangle^{s}]\big]u\big\|_{L^2}^2
 =\|\langle B(v)\xi \rangle^{-s}c(v,\xi)u\|_{L^2}^2 \\
\lesssim \|\langle v \rangle^{\gamma+2}u\|_{L^2}^2,
\end{multline}
since
$$|\langle B(v)\xi \rangle^{-s}c(v,\xi)| \lesssim \langle v \rangle^{\gamma+2},$$
uniformly with respect to the parameter $\xi$ in $\rr^3$.
It therefore follows from (\ref{ed-11}), (\ref{ed4}) and (\ref{ed7}) that
\begin{multline}\label{ed6}
\|\langle B(v)\xi \rangle^{-s}Q.\mu(v)[Q,\langle B(v)\xi \rangle^{s}]u\|_{L^2}^2 \lesssim
\|\langle B(v)\xi \rangle^{-s}b(v,\xi).Qu\|_{L^2}^2 \\ +
\big\|\langle B(v)\xi \rangle^{-s}\big[Q,\mu(v)[Q,\langle B(v)\xi \rangle^{s}]\big]u\big\|_{L^2}^2 \lesssim
\|\langle v \rangle^{\gamma+1}Qu\|_{L^2}^2+\|\langle v \rangle^{\gamma+2}u\|_{L^2}^2,
\end{multline}
since
$$|\langle B(v)\xi \rangle^{-s}b(v,\xi)| \lesssim \langle v \rangle^{\gamma+1},$$
uniformly with respect to the parameter $\xi$ in $\rr^3$.
One can then deduce from (\ref{dio1b}), (\ref{ed3}) and (\ref{ed6}) the result of Lemma~\ref{lemme2}.~$\Box$

\bigskip

\begin{lemma}\label{lemme3}
Let $P$ be the linear Landau-type operator defined in \emph{(\ref{fili11})}. Then
$$\|\langle v \rangle^{\gamma+2}u\|_{L^2}^2+\|\langle v \rangle^{\gamma+1} \nabla_vu\|_{L^2}^2+\|\langle v \rangle^{\gamma+1} Qu\|_{L^2}^2 \\ \lesssim \|Pu\|_{L^2}^2+\|u\|_{L^2}^2,$$ where $\|\cdot\|_{L^2}$ stands for the $L^2(\rr_{x,v}^{6})$-norm.
\end{lemma}

\bigskip

\noindent
\textit{Proof of Lemma~\ref{lemme3}}.
We may write that
\begin{multline}\label{j1}
\textrm{Re}(Pu,\langle v \rangle^{\gamma+2}u)=(D_v.\lambda(v)D_{v}u,\langle v \rangle^{\gamma+2}u)+(Q.\mu(v)Qu,\langle v \rangle^{\gamma+2}u) \\+\big\|\sqrt{F(v)}\langle v \rangle^{\frac{\gamma}{2}+1}u\big\|_{L^2}^2.
\end{multline}
Recalling that
$$Q=(v \wedge \eta)^w=\left(
\begin{array}{c}
v_2D_{v_3}-v_3D_{v_2}\\
v_3D_{v_1}-v_1D_{v_3}\\
v_1D_{v_2}-v_2D_{v_1}
\end{array}
\right),$$
let us notice the key commutation
$$[Q,\langle v \rangle^{\gamma+2}]=0,$$
coming from the direct computations
$$[D_{v_j},\langle v \rangle^{\gamma+2}]=\frac{\gamma+2}{i}\langle v \rangle^{\gamma}v_j,$$
when $j=1,2,3$.
It follows from (\ref{w2Q}) that
\begin{equation}\label{j2}
(Q.\mu(v)Qu,\langle v \rangle^{\gamma+2}u)=(\mu(v)Qu,\langle v \rangle^{\gamma+2}Qu) \gtrsim \|\langle v \rangle^{\gamma+1} Qu\|_{L^2}^2
\end{equation}
and
\begin{equation}\label{j3}
(\lambda(v)D_{v}u,\langle v \rangle^{\gamma+2}D_vu) \gtrsim \|\langle v \rangle^{\gamma+1} \nabla_vu\|_{L^2}^2.
\end{equation}
By writing that
\begin{multline}\label{j4}
(D_v.\lambda(v)D_{v}u,\langle v \rangle^{\gamma+2}u)=(\lambda(v)D_{v}u,[D_v,\langle v \rangle^{\gamma+2}]u)\\ +
(\lambda(v)D_{v}u,\langle v \rangle^{\gamma+2}D_vu)
\end{multline}
and noticing from (\ref{w2Q}) that
$$\|\langle v \rangle^{\gamma+2}u\|_{L^2}^2 \lesssim \big\|\sqrt{F(v)}\langle v \rangle^{\frac{\gamma}{2}+1}u\big\|_{L^2}^2 ,$$
we deduce from the Cauchy-Schwarz inequality, (\ref{eq-3.1Q}), (\ref{j1}), (\ref{j2}), (\ref{j3}) and (\ref{j4}) that
\begin{align*}
& \ \|\langle v \rangle^{\gamma+2}u\|_{L^2}^2+\|\langle v \rangle^{\gamma+1} \nabla_vu\|_{L^2}^2+\|\langle v \rangle^{\gamma+1} Qu\|_{L^2}^2 \\
\lesssim & \ \|Pu\|_{L^2}^2+
\big|\big(\lambda(v)D_{v}u,[D_v,\langle v \rangle^{\gamma+2}]u\big)\big|\\
 \lesssim & \ \|Pu\|_{L^2}^2+
\|\langle v \rangle ^{-\gamma/2}\lambda(v)D_{v}u\|_{L^2} \|\langle v \rangle ^{\gamma/2}[D_v,\langle v \rangle^{\gamma+2}]u\|_{L^2}\\
\lesssim & \ \|Pu\|_{L^2}^2+
\frac{1}{\delta}\|\langle v \rangle ^{\gamma/2}\nabla_{v}u\|_{L^2}^2+
\delta \|\langle v \rangle ^{\gamma/2}[D_v,\langle v \rangle^{\gamma+2}]u\|_{L^2}^2,
\end{align*}
for any constant $0<\delta \leq 1$. Symbolic calculus shows that
$$[D_v,\langle v \rangle^{\gamma+2}]=a(v),$$
where $a$ is a smooth function depending only on the variable $v$ and verifying
$$|a(v)| \lesssim \langle v \rangle^{\gamma+1}.$$
This implies that
$$\|\langle v \rangle ^{\gamma/2}[D_v,\langle v \rangle^{\gamma+2}]u\|_{L^2}^2 \lesssim \|\langle v \rangle^{\gamma+2}u\|_{L^2}^2,$$
since $\gamma \in [-3,1]$. By choosing the positive constant $0 <\delta \ll 1$ sufficiently small and using (\ref{and2}) to estimate from above the term
$$\|\langle v \rangle ^{\gamma/2}\nabla_{v}u\|_{L^2}^2,$$
we obtain the estimate
\begin{equation}\label{and3}
\|\langle v \rangle^{\gamma+2}u\|_{L^2}^2+\|\langle v \rangle^{\gamma+1} \nabla_vu\|_{L^2}^2+\|\langle v \rangle^{\gamma+1} Qu\|_{L^2}^2 \\ \lesssim \|Pu\|_{L^2}^2+\|u\|_{L^2}^2,
\end{equation}
which proves Lemma~\ref{lemme3}.~$\Box$

\bigskip \bigskip

\noindent \it End the proof of Proposition~\ref{propla}. \rm \ \  Working on the Fourier side in the $x$ variable, we deduce from Proposition~\ref{prop1} and Lemma~\ref{lem8} that
$$\|\langle B(v)\xi \rangle^{1/3}u\|_{L^2}^2 \lesssim \|\langle B(v)\xi \rangle^{-1/3}Pu\|_{L^2}
\|\langle B(v)\xi \rangle^{1/3}u\|_{L^2}+\|u\|_{L^2}^2,$$
uniformly with respect to the parameter $\xi \in \rr^n$; where $\|\cdot\|_{L^2}$ stands for the $L^2(\rr_{v}^{3})$-norm.
By substituting $\langle B(v)\xi \rangle^{1/3}u$ to $u$ in this estimate, we obtain that
\begin{align}\label{vienne1}
\|\langle B(v)\xi \rangle^{2/3}u\|_{L^2}^2  \lesssim & \  \|\langle B(v)\xi \rangle^{-1/3}P\langle B(v)\xi \rangle^{1/3}u\|_{L^2}^2+\|\langle B(v)\xi \rangle^{1/3}u\|_{L^2}^2\\
\lesssim & \ \|Pu\|_{L^2}^2+\|\langle B(v)\xi \rangle^{-1/3} [P,\langle B(v)\xi \rangle^{1/3}]u\|_{L^2}^2+\|u\|_{L^2}^2,
\end{align}
uniformly with respect to the parameter $\xi \in \rr^n$; since from (\ref{li4b}) and Lemma~\ref{lem8},
$$\|\langle B(v)\xi \rangle^{1/3}u\|_{L^2}^2 \lesssim \|Pu\|_{L^2}^2+\|u\|_{L^2}^2.$$
We notice from Lemma~\ref{lemme1} and Lemma~\ref{lemme2} that
\begin{multline*}
 \|\langle B(v)\xi \rangle^{-1/3}[P,\langle B(v)\xi \rangle^{1/3}]u\|_{L^2}^2
\lesssim
\|\langle B(v)\xi \rangle^{-1/3}[D_v.\lambda(v)D_v,\langle B(v)\xi \rangle^{1/3}]u\|_{L^2}^2
 +\\ \|\langle B(v)\xi \rangle^{-1/3}[Q.\mu(v)Q,\langle B(v)\xi \rangle^{1/3}]u\|_{L^2}^2 \\
 \lesssim  \|\langle v \rangle^{\gamma+1}Qu\|_{L^2}^2+\|\langle v \rangle^{\gamma+2}u\|_{L^2}^2+
 \|\langle v \rangle^{\gamma}\nabla_vu\|_{L^2}^2,
\end{multline*}
uniformly with respect to the parameter $\xi \in \rr^n$.
According to Lemma~\ref{lemme3}, by coming back to the direct side in the  $x$ variable and integrating with respect to this variable, this implies that
\begin{equation}\label{ed10}
\|\langle B(v)D_x \rangle^{-1/3}[P,\langle B(v)D_x \rangle^{1/3}]u\|_{L^2}^2 \lesssim \|Pu\|_{L^2}^2+\|u\|_{L^2}^2,
\end{equation}
where $\|\cdot\|$ stands for the $L^2(\rr_{x,v}^{6})$-norm.
We finally conclude from (\ref{w2Q}), (\ref{li4}), (\ref{vienne1}) and (\ref{ed10}) that there exists $C>0$ such that for all $u \in \mathcal{S}(\rr_{x,v}^{2n})$,
$$\|\langle B(v)D_x \rangle^{2/3}u\|_{L^2}^2  \leq C(\|Pu\|_{L^2}^2+\|u\|_{L^2}^2)$$
and
$$\|\langle v \rangle^{\gamma/3}|D_x|^{2/3}u\|_{L^2}^2 +\|\langle v \rangle^{\gamma/3}|v \wedge D_x|^{2/3}u\|_{L^2}^2  \leq C(
 \|Pu\|_{L^2}^2+\|u\|_{L^2}^2),$$
 where $\|\cdot\|_{L^2}$ stands for the $L^2(\rr_{x,v}^{6})$-norm.  This ends the proof of Proposition~\ref{propla}. \fin

\subsubsection{\underline{Velocity estimates}}

In this subsection, we begin by proving the following estimate:

\bigskip

\begin{lemma} \label{Bvxi}
Let $P$ be the linear Landau-type operator defined in \emph{(\ref{fili11})}. Then, there exists $C>0$ such that for all $u \in \mathcal{S}(\rr_{x,v}^{6})$,
$$|\Re(B(v)D_v u , B(v) D_x u)|  \leq C(
 \|Pu\|_{L^2}^2+\|u\|_{L^2}^2 + \|\langle B(v)D_v\rangle^2 u\|_{L^2}\|\langle B(v)D_x \rangle^{2/3}u\|_{L^2}),$$
where $\|\cdot\|$ stands for the $L^2(\rr_{x,v}^{6})$-norm and $\langle B(v)D_v\rangle^2 $ stands for the operator 
$$1+ D_v. B(v)^T B(v) D_v.$$
\end{lemma}

\bigskip

\proof
Proposition~\ref{propla} shows that there exists a positive constant $C>0$ such that for all $u \in \mathcal{S}(\rr_{x,v}^6)$,
\begin{equation} \label{debut}
\begin{split}
& \ |\Re(B(v)D_v u , B(v)D_x  u)| \\
& \ = \big|\Re\big(\seq{ B(v) D_x }^{1/3}  B(v)D_v u ,\seq{B(v) D_x }^{-1/3}   B(v) D_x  u\big)\big| \\
& \ \leq \|\seq{B(v) D_x }^{1/3} B(v)D_v u\|_{L^2}^2 + \|\seq{B(v) D_x }^{2/3} u\|_{L^2}^2 \\
& \ \leq ( \seq{ B(v) D_x }^{2/3}  B(v)D_v u, B(v)D_v u ) +
 C \|Pu\|_{L^2}^2+ C \|u\|_{L^2}^2,
\end{split}
\end{equation}
where $\|\cdot\|$ stands for the $L^2(\rr_{x,v}^{6})$-norm. In order to prove Lemma~\ref{Bvxi}, it only remains to estimate from above the term
$$(\seq{ B(v) D_x }^{2/3} B(v)D_v u, B(v)D_v u).$$ We work from now
in $L^2(\rr_{v}^{3})$ by considering the Fourier dual variable $\xi$ of the space variable $x$ as a parameter. We first write
\begin{align*}
& \  \big(\langle B(v)\xi\rangle^{2/3}  B(v)D_v u, B(v)D_v u\big) \\
& \ = \Re\big([\seq{ B(v) \xi }^{2/3} ,  B(v)D_v ]  u, B(v)D_v u\big) + \Re\big(  B(v)D_v \seq{ B(v) \xi }^{2/3}  u, B(v)D_v u\big) \\
& \ = \textrm{I} + \textrm{II}
\end{align*}
Let us first deal with the term II. We write
\begin{align*} %\label{part1}
\textrm{II} & \ =  \Re\big(  B(v)D_v \seq{ B(v) \xi }^{2/3}  u, B(v)D_v u\big)  =  \Re\big( \seq{ B(v) \xi }^{2/3}  u, D_v. B(v)^T  B(v)D_v u\big) \notag \\
& \  \leq \Re\big( \seq{ B(v) \xi }^{2/3}  u,(1+ D_v. B^T(v)  B(v)D_v) u\big) \notag \\
& \  = \Re\big( \seq{ B(v) \xi }^{2/3}  u,\seq{  B(v)D_v}^2 u\big) \notag \\
& \  \hspace{1cm} \leq \| \seq{ B(v) \xi }^{2/3}  u \|_{L^2} \| \seq{  B(v)D_v}^2 u \|_{L^2}.
\end{align*}
Let us now deal with the term I. For all
$j \in \set{1,2,3}$, we have
$$
(B(v)D_v)_j = \sum_{k=1}^3 B_{j,k}(v)D_{v_k} = \Big(\sum_{k=1}^3 B_{j,k}(v) \eta_k\Big)^w + i r_j(v)
$$
where $r$ is the vectorial multiplication operator with real-valued entries 
\begin{equation}\label{carm1}
r_j(v) = \frac{1}{2}\sum_{k=1}^3 (\D_k B_{j,k})(v).
\end{equation}
 With these notations, we can write
\begin{align*}
\textrm{I} & \ = \Re\big([\seq{ B(v) \xi }^{2/3} ,  (B(v) \eta)^w + ir(v) ]  u,
(B(v) \eta)^wu + i r(v) u\big).
\end{align*}
Now since $r(v)$ and  $\langle B(v)\xi\rangle^{2/3}$ are multiplication operators, they commute. Recall  the well known identity $\Re ([D,E]u, Fu) = \frac{1}{2}\Re ([F, [D,E]]u,u)$ valid when
$u \in \sss$ for the formally selfadjoint operators $E$, $F$ and $D$. We apply it with $E=F=(B(v) \eta)^w$ and $D=\seq{ B(v) \xi }^{2/3}$. It follows that
\begin{equation} \label{commm}
\begin{split}
\textrm{I}  &   = \frac{1}{2} \sum_{j=1}^3\Re\big([ (B(v) \eta)_j^w ,  [\seq{ B(v) \xi }^{2/3} ,  (B(v) \eta)_j^w ]  u,
 u\big) \\  & \ \ \ \ \ \ \ \ \ \ \ \ \ \ \ \ \ \ -   \Re\big( i r(v). [\seq{ B(v) \xi }^{2/3} ,  (B(v) \eta)^w ]u, u \big) \\
 &   \defegal  ( c_1^w u,u) +  ( c_2^w u,u).
\end{split}
\end{equation}
We shall then study each commutator appearing in the previous formula. 
By using that the Weyl symbol of the commutator $[\seq{ B(v) \xi }^{2/3} ,  (B(v) \eta)^w   ]$ is exactly given by $iB(v) \nabla_v( \seq{ B(v) \xi }^{2/3})$, another use of symbolic calculus shows that
$$c_1 =\frac{1}{2}\sum_{j=1}^3  \frac{\partial}{\partial \eta} \big((B(v) \eta)_j\big) \cdot \frac{\partial}{\partial v} \big( \big[B(v) \nabla_v( \seq{ B(v) \xi }^{2/3}) \big]_j\big).$$
This is a multiplication operator. Lemma~\ref{lem0} together with (\ref{eq-1QQ}) show that 
\begin{align} \label{cun}
|c_1|  \lesssim |B(v)|^2 \seq{ B(v) \xi }^{2/3} + |B(v)||B'(v)|\seq{ B(v) \xi }^{2/3} \lesssim \langle v \rangle^{\gamma+2}\seq{ B(v) \xi }^{2/3},
\end{align}
uniformly with respect to the parameter $\xi \in \rr^3$. 
As a consequence, it follows that 
$$
 |( c_1^w u,u)| \lesssim   ( \seq{ B(v) \xi }^{2/3}u,\seq{v}^{\gamma+2} u)_{L^2} \lesssim  \|\seq{ B(v) \xi }^{2/3}u\|_{L^2} \|\seq{v}^{\gamma+2} u\|_{L^2}.
$$
Notice now that 
$$c_2(v)=r(v).B(v) \nabla_v( \seq{ B(v) \xi }^{2/3}),$$
since $r$ and $B$ have real-valued entries. This is again a multiplication operator whose symbol may be bounded from above as
\begin{align*}  %\label{cdeux}
|c_2|  \lesssim |B(v)||B'(v)|\seq{ B(v) \xi }^{2/3}  \lesssim  \seq{v}^{\gamma +2}\seq{ B(v) \xi }^{2/3},
\end{align*}
uniformly with respect to the parameter $\xi \in \rr^3$; according to Lemma~\ref{lem0}, (\ref{eq-1QQ}) and (\ref{carm1}). Proceeding as for $c_1$, we obtain the second estimate 
$$|( c_2^w u,u)| \lesssim  \|\seq{ B(v) \xi }^{2/3}u\|_{L^2} \|\seq{v}^{\gamma+2} u\|_{L^2},$$
which implies that 
\begin{equation*} %\label{II}
|\textrm{I}| \lesssim \|\seq{ B(v) \xi }^{2/3}u\|_{L^2} \|\seq{v}^{\gamma+2} u\|_{L^2}.
\end{equation*}
According to the estimates of the two terms I and II,  we obtain 
after integration in the $\xi$ side the new estimate in $L^2(\R^6_{x,v})$,
 \begin{align*}
& \ (\seq{ B(v) D_x }^{2/3}  B(v)D_v u, B(v)D_v u) \\
& \ \lesssim  \|\seq{ B(v) D_x }^{2/3}u\|_{L^2} \|\seq{v}^{\gamma+2} u\|_{L^2}+\|\seq{ B(v) D_x}^{2/3}u\|_{L^2} \|\seq{B(v) D_v}^2 u\|_{L^2}.
\end{align*}
%Putting this in (\ref{debut})
Proposition~\ref{propla} and Lemma~\ref{lemme3} give
 \begin{align*}
& \ (\seq{ B(v) D_x }^{2/3}  B(v)D_v u, B(v)D_v u) \\
& \ \lesssim  \ \|\seq{ B(v) D_x }^{2/3}u\|_{L^2}^2+  \|\seq{v}^{\gamma+2} u\|_{L^2}^2+\|\seq{ B(v) D_x }^{2/3}u\|_{L^2} \|\seq{B(v) D_v}^2 u\|_{L^2} \\
& \  \lesssim   \|Pu\|_{L^2}^2+\|u\|_{L^2}^2+\|\seq{ B(v) D_x }^{2/3}u\|_{L^2} \|\seq{B(v) D_v}^2 u\|_{L^2},
\end{align*}
which together with (\ref{debut}) finally complete the proof of this lemma.  \fin

\bigskip

\noindent
We now prove a result fully independent of the $x$ variable. The following proof relies on the use of the Fefferman-Phong inequality.

\bigskip

 \begin{lemma} \label{feff}
 Let $B(v)$ be the matrix defined in \emph{(\ref{an2QQ})} and denote again 
 $$\langle B(v)D_v\rangle^2 = 1+ D_v. B(v)^T B(v) D_v.$$ 
 Then, there exists $C >0$ such that for all $u \in \sss(\R^3_v)$,
 $$
 \|D_v. \seq{v}^\gamma D_v u\|_{L^2}^2 \leq C \|\langle B(v)D_v\rangle^2 u\|_{L^2}^2 + C
\|\seq{v}^{\gamma+1} D_v u\|^2_{L^2} + C\|\seq{v}^{\gamma+1}  u\|^2_{L^2},$$
where $\|\cdot\|_{L^2}$ stands for the $L^2(\rr_v^3)$-norm.
 \end{lemma}

\bigskip

 \proof
 Recalling from (\ref{eq-2QQ}) that $A(v) = B(v)^T B(v)$, we may rewrite the terms 
 $$\|D_v. \seq{v}^\gamma D_v u\|_{L^2}^2=\sep{ ( D_v. \seq{v}^\gamma D_v)^2 u , u }$$
 and 
 $$ \|\langle B(v)D_v\rangle^2 u\|_{L^2}^2= \sep{ (1+ D_v. A(v) D_v)^2 u, u}.$$
  %estimate to be proved as
 %\begin{equation} \label{vv}
 %\sep{ ( D_v. \seq{v}^\gamma D_v)^2 u , u } \leq C \sep{ (\textrm{Id}+ D_v. A(v) D_v)^2 u, u} + C \norm{u}^2.
 %\end{equation}
We introduce the following metric
$$
\widetilde{\Gamma} = \frac{dv^2}{\langle v \rangle^2} + \frac{d\eta^2}{\langle \eta \rangle^2}.
$$
It is easy to check that this metric $\tilde{\Gamma}$ is admissible (slowly varying, satisfying the uncertainty principle and temperate) with gain
$$\lambda_{\tilde{\Gamma}}(v,\eta)=\langle v \rangle \langle \eta \rangle.$$ 
Let $a(v,\eta)$, respectively $\tilde{a}(v,\eta)$; be the Weyl symbol of the operator $(D_v. \seq{v}^\gamma D_v)^2$, respectively the Weyl symbol of the operator $D_v. \seq{v}^\gamma D_v$. Notice that 
$$a \in  S( \seq{\eta}^4 \seq{v}^{2 \gamma}, \widetilde{\Gamma}) \textrm{ and } \tilde{a} \in  S( \seq{\eta}^2 \seq{v}^{ \gamma}, \widetilde{\Gamma}).$$
Symbolic calculus shows
$$a_1 \defegal a - \tilde{a}^2 -\frac{1}{2i}\{\tilde{a},\tilde{a}\}=a - \tilde{a}^2  \in S( \seq{\eta}^2\seq{v}^{2\gamma-2}, \widetilde{\Gamma}).$$
It follows that
$$a(v, \eta) - a_1(v,\eta) \lesssim \langle \eta \rangle^4 \seq{v}^{2 \gamma}.$$
Let $b(v,\eta)$, respectively $\tilde{b}(v,\eta)$;  be the Weyl symbol of the operator 
$$(1+D_v. A(v) D_v)^2,$$ 
respectively the Weyl symbol of the operator $1+D_v. A(v) D_v$. Notice from (\ref{eq-3QQ}) that
$$b \in  S( \seq{\eta}^4 \seq{v}^{2 \gamma+4}, \widetilde{\Gamma}) \textrm{ and } \tilde{b} \in  S( \seq{\eta}^2 \seq{v}^{ \gamma+2}, \widetilde{\Gamma}).$$
Symbolic calculus shows 
$$b_1\defegal b -\tilde{b}^2-\frac{1}{2i}\{\tilde{b},\tilde{b}\}=b - \tilde{b}^2  \in S( \seq{\eta}^2\seq{v}^{2\gamma+2}, \widetilde{\Gamma}).$$
A direct computation using symbolic calculus, (\ref{eq-3QQ}) and (\ref{eq0QQ}) shows that there exists a positive constant $C>0$ such that
\begin{multline*}
\tilde{b}(v,\eta)=1+\sum_{j,k=1}^3\Big[ A_{j,k}(v)\eta_j\eta_k+\frac{i}{2}\eta_j \partial_k A_{j,k}(v)-\frac{i}{2}\eta_k\partial_jA_{j,k}(v)+\frac{1}{4}\partial_{j,k}^2A_{j,k}(v) \Big] \\
=1+\sum_{j,k=1}^3\Big[ A_{j,k}(v)\eta_j\eta_k+\frac{1}{4}\partial_{j,k}^2A_{j,k}(v)\Big] \geq c \langle v \rangle^{\gamma}|\eta|^2-C\langle v \rangle^{\gamma},
\end{multline*}
since by symmetry $A_{j,k}=A_{k,j}$, for any $1 \leq j,k \leq 3$.
It follows that one can find a new positive constant $C>0$ such that
$$b(v,\eta)  - b_1(v, \eta)=\tilde{b}(v,\eta)^2 \gtrsim   \abs{\eta}^4 \seq{v}^{2 \gamma} - C \langle v \rangle^{2\gamma}.$$
The Feffermann-Phong inequality therefore yields
\begin{equation} \label{f1}
(a^w u, u) \lesssim (b^w u,u) + ((a_1-b_1)^w u,u) + (r^w u, u)
\end{equation}
with $r$ and $a_1-b_1 \in S(\seq{\eta}^2 \seq{v}^{2\gamma + 2}, \widetilde{\Gamma})$.
We iterate this method in order to treat the symbols $a_1-b_1$ and $r$. Define
$$c(v,\eta) \defegal \sigma( D_v. \seq{v}^{2\gamma + 2} D_v) + \seq{v}^{2\gamma + 2},$$
where $\sigma( D_v. \seq{v}^{2\gamma + 2} D_v)$ stands for the Weyl symbol of the operator $D_v. \seq{v}^{2\gamma + 2} D_v$.
A new straightforward computation shows that
$$c(v,\eta) = \seq{v}^{2\gamma + 2} \seq{ \eta}^2  + c_1(v),$$
with $|c_1(v) | \lesssim \seq{v}^{2\gamma}$.
Notice that
$$c(v, \eta) \in S( \seq{\eta}^2 \seq{v}^{2 \gamma+2}, \widetilde{\Gamma}) \ \ \ \textrm{ and}
\ \ \ \
c(v,\eta)   \gtrsim   \seq{\eta}^2 \seq{v}^{2 \gamma+2} - C.$$
By using the Fefferman-Phong inequality, we obtain that
\begin{equation} \label{f2}
((a_1-b_1+r)^w u,u) \lesssim (c^w u, u) + (s^w u, u)+\|u\|_{L^2}^2,
\end{equation}
with $s \in S( \seq{v}^{2\gamma}, \widetilde{\Gamma})$. We make a last iteration of the previous analysis and define
$$d(v) = \seq{v}^{2\gamma}.$$
By using again the Fefferman-Phong inequality, we obtain that there exists a new positive constant $C>0$ such that
\begin{equation} \label{f3}
(s^w u,u) \leq (d^w u, u) + C\|u\|^2_{L^2},
\end{equation}
since $\gamma \in [-3,1]$. 
Putting all together estimates (\ref{f1}), (\ref{f2}) and (\ref{f3}) provides the estimate
\begin{align*}
 & \quad \big((D_v. \seq{v}^\gamma D_v)^2 u , u\big) \\
 & \ \lesssim (b^w u,u) + (c^w u,u)  + (d^w u, u) +  \|u\|^2_{L^2} \\
 & \ \lesssim \big((1 + D_v. A(v) D_v)^2 u, u\big) + \|\seq{v}^{\gamma+1} D_v u\|^2_{L^2} + \|\seq{v}^{\gamma+1}  u\|^2_{L^2} + \|\seq{v}^{\gamma}  u\|^2_{L^2} + \|u\|^2_{L^2} \\
 & \ \lesssim \big(\seq{ B(v) D_v}^4 u, u\big) +\|\seq{v}^{\gamma+1} D_v u\|^2_{L^2} + \|\seq{v}^{\gamma+1}  u\|^2_{L^2},
\end{align*}
which proves this lemma.  \fin

\bigskip

\begin{proposition}\label{propla2}
Let $P$ be the linear Landau-type operator defined in \emph{(\ref{fili11})}. Then, there exists $C>0$ such that for all $u \in \mathcal{S}(\rr_{x,v}^{6})$,
$$\|\langle v \rangle^{\gamma}|D_v|^{2}u\|_{L^2}^2+\|\langle v \rangle^{{\gamma}}|v \wedge D_v|^{2}u\|_{L^2}^2  \leq C(
 \|Pu\|_{L^2}^2+\|u\|_{L^2}^2),$$
where $\|\cdot\|_{L^2}$ stands for the $L^2(\rr_{x,v}^{6})$-norm.
\end{proposition}

\bigskip

\proof
As a first step, we shall prove the following estimate
\begin{equation} \label{step1}
 \|\seq{ B(v)  D_v}^2 u\|_{L^2}^2  \leq C(
 \|Pu\|_{L^2}^2+\|u\|_{L^2}^2),
 \end{equation}
 with $\langle B(v)D_v\rangle^2 = 1+ D_v. B(v)^T B(v) D_v$. Recalling that 
 $$P=iv.D_x+D_v.B(v)^TB(v)D_v+F(v),$$
 we may write for any $u \in \mathcal{S}(\rr_{x,v}^{6})$,
\begin{equation*}
\begin{split}
  \big(\seq{ B(v)  D_v}^4 u,  u\big) & = \big(\seq{B(v)  D_v}^2 u,  u\big) + \big(\seq{B(v)  D_v}^2 D_v . B(v)^T B(v) D_v u, u\big) \\
  &   \leq \Re(Pu,u) + \|u\|_{L^2}^2 +\Re(\seq{B(v)  D_v}^2 (P - F(v) - iv.D_x) u, u) \\
&   \leq \Re(Pu,u) + \|u\|_{L^2}^2 + \frac{1}{4} \|\seq{B(v)  D_v}^2 u\|_{L^2}^2  \\
& \ \ \ \ \ \ \  +2 \|P u\|_{L^2}^2 + 2 \|F(v)u\|_{L^2}^2
   - \Re(\seq{B(v)  D_v}^2 u, iv.D_x u). \\
\end{split}
\end{equation*}
By using Lemma \ref{lemme3} and (\ref{eq-3.1Q}), we obtain that
\begin{equation} \label{cool}
\begin{split}
  (\seq{ B(v)  D_v}^4 u,  u) & \leq C \big(\|P u\|_{L^2}^2 + \|u\|_{L^2}^2 +
   |\Re(\seq{B(v) D_v}^2 u, iv.D_x u)|\big).
\end{split}
\end{equation}
Noticing that the operator $iv.D_x$ is formally skew-adjoint on $L^2$, a direct computation gives that
\begin{align*}
& \Re(\seq{B(v)  D_v}^2 u, iv.D_x u) \\
& \ \ \ \ \ \ =  \  \frac{1}{2} \Re([\seq{B(v)  D_v}^2, iv.D_x]u,u)\\
& \ \ \ \ \ \ =  \ \frac{1}{2} \Re([ D_v .B(v)^T B(v)  D_v , i v. D_x] u , u) \\
& \ \ \ \ \ \ =  \ \frac{1}{2} \Re\sep{ ([ D_v, i v. D_x]. B(v)^T B(v)  D_v + D_v. B(v)^T B(v) [D_v, iv. D_x] )u , u} \\
& \ \ \ \ \ \ =  \ \frac{1}{2} \Re \sep{ ( D_x. B(v)^T B(v)  D_v + D_v. B(v)^T B(v)  D_x ) u , u} \\
& \ \ \ \ \ \  =  \ \Re(B(v)  D_v u , B(v)D_x u).
\end{align*}
It then follows from (\ref{cool}) and Lemma~\ref{Bvxi} that there exists a new positive constant $C>0$ such that for all $u \in \mathcal{S}(\rr_{x,v}^{6})$,
 \begin{equation*}
\begin{split}
  & \quad (\seq{ B(v)  D_v}^4u,u) \\
  & \lesssim \|P u\|_{L^2}^2 + \|u\|_{L^2}^2 +|\Re(\seq{B(v)  D_v}^2 u, iv.D_x u)| \\
  &   \leq C\big(\|P u\|_{L^2}^2 + \|u\|_{L^2}^2+\|\seq{ B(v)D_x}^{2/3}u\|_{L^2}^2\big) +
  \frac{1}{2}  \|\seq{ B(v)  D_v}^2 u\|_{L^2}^2
\end{split}
\end{equation*}
By using Proposition \ref{propla} to estimate from above the third term, we obtain that one can find a new positive constant $C>0$ such that for all $u \in \mathcal{S}(\rr_{x,v}^{6})$,
$$\|\seq{ B(v)  D_v}^2 u\|_{L^2} \leq C(\|P u\|_{L^2}^2 + \|u\|_{L^2}^2) +  \frac{1}{2}  \|\seq{ B(v)  D_v}^2 u\|_{L^2}.$$
This proves (\ref{step1}). We now deal with the core of the proof of Proposition~\ref{propla2}. We first write that
\begin{align*}
& \|\seq{v}^\gamma |D_v|^2 u\|^2_{L^2} \leq 2 \|D_v.\seq{v}^\gamma D_v u\|^2_{L^2} + 2
 \|[\seq{v}^\gamma, D_v]. D_v u\|^2_{L^2} \\
 & \ \leq 2 \|D_v.\seq{v}^\gamma D_v u\|^2_{L^2} + C  \|\seq{v}^{\gamma-1} D_v u\|^2_{L^2} \\
  & \ \lesssim \|\seq{ B(v) D_v}^2 u\|^2_{L^2} + \|\seq{v}^{\gamma+1} D_v u\|^2_{L^2} + \|\seq{v}^{\gamma+1}  u\|^2_{L^2} + \|u\|^2_{L^2}
\end{align*}
where we used  $\seq{v}^{\gamma-1} \leq \seq{v}^{\gamma+1}$  and
Lemma \ref{feff} in a crucial way.
Using then  inequality (\ref{step1}) and Lemma \ref{lemme3}, we get the following  result
\begin{equation} \label{step2}
\|\seq{v}^\gamma |D_v|^2 u\|^2_{L^2} \lesssim \|P u\|_{L^2}^2 + \|u\|_{L^2}^2.
\end{equation}
We now deal with the term with a cross product. Recalling (\ref{eq-3.1Q}), (\ref{w2Q}) and (\ref{ed0}), we may write
\begin{equation} \label{trucc}
\begin{split}
& \|\seq{v}^\gamma |v \wedge D_v|^2 u\|^2_{L^2} \leq \|\mu(v) |v \wedge D_v|^2 u\|^2_{L^2} \\
& \ \leq 2\|(v \wedge D_v) .\mu(v) (v \wedge D_v) u\|^2_{L^2} + 2\|[\mu(v), v \wedge D_v] . (v \wedge D_v) u\|^2_{L^2} \\
& \ \leq 2\|(v \wedge D_v) .\mu(v) (v \wedge D_v) u\|^2_{L^2} +  C \|\seq{v}^{\gamma+1} D_v u\|^2_{L^2}
\end{split}
\end{equation}
where we used  that $[\mu(v), v \wedge D_v] . (v \wedge D_v) = R(v) D_v$; with $R(v)$ an explicit matrix whose entries are all bounded by a positive constant times the function $\seq{v}^{\gamma +1}$. Recall that
$$\seq{ B(v) D_v}^2 =1 +  D_v. B(v)^T B(v) D_v = 1 +D_v.\lambda(v)D_{v}+(v \wedge D_v).\mu(v) (v \wedge D_v).
$$
As a consequence, we deduce from (\ref{trucc}) that
\begin{multline*}
\|\seq{v}^\gamma |v \wedge D_v|^2 u\|^2_{L^2}  \lesssim \|\seq{ B(v) D_v}^2 u\|^2_{L^2} +
\|D_v.\lambda(v) D_v u\|^2_{L^2} + \|u\|^2_{L^2} \\ +  \|\seq{v}^{\gamma+1} D_v u\|^2_{L^2}.
\end{multline*}
Another direct argument of commutation using (\ref{eq-3.1Q}) gives the estimate
$$ \|D_v.\lambda(v) D_v u\|^2_{L^2} \lesssim \|\seq{v}^\gamma |D_v|^2 u\|^2_{L^2} + \|\seq{v}^{\gamma-1}  D_vu\|^2_{L^2},$$
which implies that  
 \begin{multline*}
 \|\seq{v}^\gamma |v \wedge D_v|^2 u\|^2_{L^2}  
   \lesssim \|\seq{ B(v) D_v}^2 u\|^2_{L^2} +
\|\seq{v}^\gamma |D_v|^2 u\|^2_{L^2} + \|u\|^2_{L^2} \\ + \|\seq{v}^{\gamma+1} D_v u\|^2_{L^2},
\end{multline*}
because $\seq{v}^{\gamma-1} \leq \seq{v}^{\gamma+1}$. One can then deduce from Lemma \ref{lemme3},  (\ref{step1}) and  (\ref{step2}) that 
\begin{align} \label{step4}
 \|\seq{v}^\gamma |v \wedge D_v|^2 u\|^2_{L^2}  \lesssim \|P u\|_{L^2}^2 + \|u\|_{L^2}^2.
\end{align}
Proposition~\ref{propla2} then directly follows from (\ref{step2}) and (\ref{step4}). \fin

\subsubsection{Proof of Theorem~\ref{mainlandau}}
Theorem~\ref{mainlandau} is now a direct consequence of Proposition~\ref{propla}, Proposition~\ref{propla2} and Lemma~\ref{lemme3}. This ends the proof of Theorem~\ref{mainlandau}.
\fin

\bigskip

\noindent
The proof of Theorem~\ref{mainlandau} can easily be adapted to obtain the following time dependent hypoelliptic estimate.

\bigskip

\begin{proposition} \label{proplat}  Let $P$ be the linear Landau-type operator defined in \emph{(\ref{fili11})}. Then, there exists a positive constant $C > 0$ such that for all $u \in \sss(\R^{7}_{t,x,v})$,
 \begin{multline*}
\|\langle v \rangle^{\gamma +2} u\|_{L^2}^2  
+  \|\seq{v}^{\gamma} |D_v|^{2} u \|_{L^2}^2
  +  \| \seq{v}^{\gamma} |v \wedge D_v|^{2} u \|_{L^2}^2 \\
+ \|\seq{v}^{\gamma/3} |D_x|^{2/3} u \|_{L^2}^2 +
  \| \seq{v}^{\gamma/3} |v \wedge D_x|^{2/3} u \|_{L^2}^2  
  \leq C (\|\D_t u + Pu\|_{L^2}^2 + \|u\|_{L^2}^2),
\end{multline*}
where $\|\cdot\|_{L^2}$ stands for the $L^2(\rr_{t,x,v}^{7})$-norm.
\end{proposition}

\bigskip

\noindent
\textit{Proof of Proposition \ref{proplat}}.
It is sufficient to notice that through all the proof of Theorem \ref{mainlandau}, one can substitute without any change the operator $\tilde{P}=i\tau+P$ to the linear Landau-type operator $P$. Indeed, the real parameter $\tau$ disappears in all the commutators involved in this analysis. Same remark when we take the real part of the $L^2$ scalar product
$$\textrm{Re}(i\tau u+Pu,u)=\textrm{Re}(Pu,u)$$
and (see (\ref{eq7})),
$$\textrm{Re}\big(i\tau u+Pu,(1-\eps G)u\big)=\textrm{Re}\big(Pu,(1-\eps G)u\big),$$
since the multiplier $G=\tilde{g}^w$, whose Weyl symbol is real-valued is a formally selfadjoint operator on $L^2$.
Proposition~\ref{proplat} then follows from the same proof as the one given for Theorem~\ref{mainlandau} after substituting the operator $\tilde{P}$ to $P$; and then coming back to the direct side in the $t$ variable; and integrating those estimates with respect to this variable. \fin

\section{Appendix on Wick calculus}\label{la3}

The purpose of this section is to recall the definition and basic properties of the Wick quantization. We follow here the presentation of the Wick quantization given by N.~Lerner in \cite{Ler03} (see also \cite{Le}); and refer the reader to his work for the proofs of the results recalled below.

The main property of the Wick quantization is its property of positivity, i.e., that non-negative Hamiltonians define non-negative operators
$$a \geq 0 \Rightarrow a^{\textrm{Wick}} \geq 0.$$
We recall that this is not the case for the Weyl quantization and refer to \cite{Ler03} for an example of non-negative Hamiltonian defining an operator which is not non-negative.

Before defining properly the Wick quantization, we first need to recall the definition of the wave packets transform of a function $u \in \mathcal{S}(\rr^n)$,
$$Wu(y,\eta)=(u,\varphi_{y,\eta})_{L^2(\rr^n)}=2^{n/4}\int_{\rr^n}{u(x)e^{- \pi (x-y)^2}e^{-2i \pi(x-y).\eta}dx}, \ (y,\eta) \in \rr^{2n}.$$
where
$$\varphi_{y,\eta}(x)=2^{n/4}e^{- \pi (x-y)^2}e^{2i \pi (x-y).\eta}, \ x \in \mathbb{R}^n;$$
and $x^2=x_1^2+...+x_n^2$. With this definition, one can check (See Lemma 2.1 in \cite{Ler03}) that
the mapping $u \mapsto Wu$ is continuous from $\mathcal{S}(\rr^n)$ to $\mathcal{S}(\rr^{2n})$, isometric from $L^{2}(\rr^n)$ to $L^2(\rr^{2n})$ and that we have the
reconstruction formula
\begin{equation}\label{lay0.1}
\forall u \in \mathcal{S}(\rr^n), \forall x \in \rr^n, \ u(x)=\int_{\rr^{2n}}{Wu(y,\eta)\varphi_{y,\eta}(x)dyd\eta}.
\end{equation}
By denoting $\Sigma_Y$ the operator defined in the Weyl quantization by the symbol
$$p_Y(X)=2^n e^{-2\pi|X-Y|^2}, \ Y=(y,\eta) \in \rr^{2n};$$
which is a rank-one orthogonal projection,
$$\big{(}\Sigma_Y u\big{)}(x)=Wu(Y)\varphi_Y(x)=(u,\varphi_Y)_{L^2(\rr^n)}\varphi_Y(x),$$
we define the Wick quantization of any $L^{\infty}(\rr^{2n})$  symbol $a$ as
\begin{equation}\label{lay0.2}
a^{\textrm{Wick}}=\int_{\rr^{2n}}{a(Y)\Sigma_Y dY}.
\end{equation}
More generally, one can extend this definition when the symbol $a$ belongs to $\mathcal{S}'(\rr^{2n})$ by defining the operator $a^{\textrm{Wick}}$ for any $u$ and $v$ in $\mathcal{S}(\rr^{n})$ by
$$<a^{\textrm{Wick}}u,\overline{v}>_{\mathcal{S}'(\rr^{n}),\mathcal{S}(\rr^{n})}=<a(Y),(\Sigma_Yu,v)_{L^2(\rr^n)}>_{\mathcal{S}'(\rr^{2n}),\mathcal{S}(\rr^{2n})},$$
where $<\textrm{\textperiodcentered},\textrm{\textperiodcentered}>_{\mathcal{S}'(\rr^n),\mathcal{S}(\rr^n)}$ denotes the duality bracket between the
spaces $\mathcal{S}'(\rr^n)$ and $\mathcal{S}(\rr^n)$. The Wick quantization is a positive quantization
\begin{equation}\label{lay0.5}
a \geq 0 \Rightarrow a^{\textrm{Wick}} \geq 0.
\end{equation}
In particular, real Hamiltonians get quantized in this quantization by formally self-adjoint operators and one has (See Proposition 3.2 in \cite{Ler03}) that $L^{\infty}(\rr^{2n})$ symbols define bounded operators on $L^2(\rr^n)$ such that
\begin{equation}\label{lay0}
\|a^{\textrm{Wick}}\|_{\mathcal{L}(L^2(\rr^n))} \leq \|a\|_{L^{\infty}(\rr^{2n})}.
\end{equation}
According to Proposition~3.3 in~\cite{Ler03}, the Wick and Weyl quantizations of a symbol $a$ are linked by the following identities
\begin{equation}\label{lay1bis}
a^{\textrm{Wick}}=\tilde{a}^w,
\end{equation}
with
\begin{equation}\label{lay2bis}
\tilde{a}(X)=\int_{\rr^{2n}}{a(X+Y)e^{-2\pi |Y|^2}2^ndY}, \ X \in \rr^{2n};
\end{equation}
and
\begin{equation}\label{lay1}
a^{\textrm{Wick}}=a^w+r(a)^w,
\end{equation}
where $r(a)$ stands for the symbol
\begin{equation}\label{lay2}
r(a)(X)=\int_0^1\int_{\rr^{2n}}{(1-\theta)a''(X+\theta Y)Y^2e^{-2\pi |Y|^2}2^ndYd\theta}, \ X \in \rr^{2n};
\end{equation}
if we use here the normalization chosen in \cite{Ler03} for the Weyl quantization
\begin{equation}\label{lay3}
(a^wu)(x)=\int_{\rr^{2n}}{e^{2i\pi(x-y).\xi}a\Big(\frac{x+y}{2},\xi\Big)u(y)dyd\xi},
\end{equation}
which differs from the one chosen in the rest of this paper.
We also recall the following composition formula obtained in the proof of Proposition~3.4 in~\cite{Ler03},
\begin{equation}\label{lay4}
a^{\textrm{Wick}} b^{\textrm{Wick}} =\Big{[}ab-\frac{1}{4 \pi} a'.b'+\frac{1}{4i \pi}\{a,b\} \Big{]}^{\textrm{Wick}}+S,
\end{equation}
with $\|S\|_{\mathcal{L}(L^2(\rr^n))} \leq d_n \|a\|_{L^{\infty}}\gamma_{2}(b),$
when $a \in L^{\infty}(\rr^{2n})$ and $b$ is a smooth symbol satisfying
$$\gamma_2(b)=\sup_{X \in \rr^{2n}, \atop T \in \rr^{2n}, |T|=1}|b^{(2)}(X)T^2| < +\infty.$$
The term $d_n$ appearing in the previous estimate stands for a positive constant depending only on the dimension $n$, and the notation $\{a,b\}$ denotes the Poisson bracket
$$\{a,b\}=\frac{\partial a}{\partial \xi}.\frac{\partial b}{\partial x}-\frac{\partial a}{\partial x}.\frac{\partial b}{\partial \xi}.$$

\bigskip
\bigskip

\noindent
\textbf{Acknowledgement.} This work was initiated during the program \textit{Selected topics in spectral theory} organized by B.~Helffer, T.~Hoffmann-Ostenhof and A.~Laptev at the Erwin Schr\"odinger Institute for Mathematical Physics, in Vienna during the summer 2009. The authors would like to thank the Institute and the organizers very much for their hospitality and the exceptional working surroundings.

\end{document}